\documentclass[times, 5p, fleqn]{elsarticle}

\usepackage{amssymb}

\usepackage{url}
\makeatletter
\AtBeginDocument{%
  \def\doi#1{\upshape\url{https://doi.org/#1}}}
\makeatother

\usepackage{tabularx}
\usepackage{mathtools}
\usepackage{amsopn}
\DeclareMathOperator{\diag}{diag}
\DeclareMathOperator{\sinc}{sinc}
\DeclareMathOperator*{\argmin}{argmin}

\DeclarePairedDelimiter{\norm}{\|}{\|}
\DeclarePairedDelimiter{\abs}{\lvert}{\rvert}
\DeclarePairedDelimiter{\lidx}{[}{]}

\graphicspath{{./}{./graphics}}

\journal{Journal of Computational Science}

\begin{document}
\begin{frontmatter}



\title{Customizable Adaptive Regularization Techniques for B-Spline Modeling\tnoteref{funding}}
\tnotetext[funding]{This work is supported by the U.S. Department of Energy, Office of Science, Advanced Scientific Computing Research under Contract DE-AC02-06CH11357, Program Manager Margaret Lentz.}

\author[argonne]{David Lenz\corref{cor}}
\ead{dlenz@anl.gov}
\author[google]{Raine Yeh}
\author[argonne]{Vijay Mahadevan}
\author[argonne]{Iulian Grindeanu}
\author[argonne]{Tom Peterka}
\cortext[cor]{Corresponding author.}

\affiliation[argonne]{organization={Mathematics and Computer Science Division, Argonne National Laboratory},
            city={Lemont},
            postcode={60439}, 
            state={IL},
            country={USA}}
\affiliation[google]{organization={Google},
            city={New York City},
            state={NY},
            country={USA}}

\begin{abstract}
B-spline models are a powerful way to represent scientific data sets with a functional approximation. However, these models can suffer from spurious oscillations when the data to be approximated are not uniformly distributed. Model regularization (i.e., smoothing) has traditionally been used to minimize these oscillations; unfortunately, it is sometimes impossible to sufficiently remove unwanted artifacts without smoothing away key features of the data set. In this article, we present a method of model regularization that preserves significant features of a data set while minimizing artificial oscillations. Our method varies the strength of a smoothing parameter throughout the domain automatically, removing artifacts in poorly-constrained regions while leaving other regions unchanged. The proposed method selectively incorporates regularization terms based on first and second derivatives to maintain model accuracy while minimizing numerical artifacts. The behavior of our method is validated on a collection of two- and three-dimensional data sets produced by scientific simulations. In addition, a key tuning parameter is highlighted and the effects of this parameter are presented in detail. This paper is an extension of our previous conference paper at the 2022 International Conference on Computational Science (ICCS)~\cite{lenz2022adaptive}. 
\end{abstract}

\end{frontmatter}

\section{Introduction}
Data sets assembled from scientific simulations or experimental readings are often defined as a list of position-value pairs, where each data point consists of a measurement and the corresponding location of that measurement. These point locations can form structured grids, unstructured meshes, or unconnected point clouds, depending on the application. Methods for analyzing these data often apply only to particular layouts; usually, numerical analysis techniques become more complex as the geometry of the point locations becomes more general (e.g. point clouds). Even seemingly straightforward tasks such as interpolation can be computationally burdensome on unstructured point clouds and numerically inaccurate on highly nonuniform meshes. One way to avoid these challenges is by representing a data set with a mathematical function and then analyzing the function instead of the original data.  This can substantially streamline the process of interpolation and differentiation away from data points, simplify visualization tasks, and make resampling the data almost trivial.

The focus of this article is the representation of scientific data sets with (tensor-product) B-splines. B-splines are a family of highly-regular functions commonly used inaccurate geometric modeling~\cite{lin2018geometric} and form the underpinnings of isogeometric analysis (IGA)~\cite{hughes2005iga}. Recent study has shown that large, complex data sets produced by scientific simulations at extreme scale can be effectively modeled by B-splines~\cite{peterka_ldav18}. Similar results have also been obtained for nonuniform rational B-splines, or NURBS, which are a generalization of B-splines~\cite{nashed2019}.

B-splines have a number of properties that make them useful as a functional representation of data. 
B-splines are high-order approximants, and evaluating, differentiating, and integrating a B-spline model is fast and numerically stable~\cite{deBoor1978guide}. Crucially, differentiation and integration can be computed in closed-form and incur no additional loss of accuracy, unlike finite differences or Riemann sums. Thus, once a sufficiently accurate spline has been computed to represent the data, it is often more productive to analyze the functional model than the original data.

However, computing a best-fit B-spline requires solving a linear system which may be ill-conditioned or rank-deficient. A common cause of this ill-conditioning is an input data set that contains both sparse and dense patches of points in proximity to each other. 
Without additional effort, solving this system can produce a function that oscillates strongly between data points or even diverges in regions where input data are very sparse. This problem can be hard to detect automatically, since error metrics are usually defined in terms of the pointwise error between the original data and the model. Spurious oscillations occurring away from the input data will not be captured by these metrics. 

To address these challenges, we developed a new method for fitting B-splines to unstructured data that reduces or eliminates oscillations while leaving critical features of the data set unchanged. Our method regularizes the solution to the B-spline fitting problem by adding a variable-strength smoothing parameter that automatically adapts based on characteristics of the input data set. This additional term smooths out spike artifacts in regions where the data set is very sparse but does not do any smoothing where data points are densely packed, thereby preserving accuracy in these regions. In addition, our method creates well-defined spline models even for data sets with irregular boundaries. No knowledge of the boundary is required; the method automatically handles areas outside the boundary that contain no data points. 

The remainder of this paper is organized as follows. A review of related ideas and methods is given in Section~\ref{sec:related-work}. In Section~\ref{sec:background}, we provide a primer on the mathematical details used to describe B-splines throughout the paper. Our main result, a method for adaptive regularization of  B-spline models, is described in Section~\ref{sec:method}. We then exhibit the performance of this method in Section~\ref{sec:results} with a series of numerical examples. We summarize directions for further research in Section~\ref{sec:future-work} and present conclusions in Section~\ref{sec:conclusions}.

\section{Related Work}\label{sec:related-work}
Creating B-spline models to represent unstructured data sets is a particular example of scattered data approximation (SDA), a broad area of study concerned with defining continuous functions that interpolate or approximate spatially scattered inputs. SDA is often applied to image reconstruction problems, where an experimental or physical constraint prohibits the collection of uniformly-spaced samples, such as medical~\cite{arigovindan2006full}, seismic~\cite{duijndam1999reconstruction}, or astronomical~\cite{vio2000reconstruction} imaging. An introductory comparison of SDA methods was compiled by Francis et al.~\cite{francis2018scattered}.

Ill-conditioned numerical methods are a persistent challenge throughout the SDA literature, and a number of techniques have been proposed to increase numerical stability. Our approach is most similar to the variational methods of SDA, in which the magnitude of the approximating function's derivative (or ``roughness'') is minimized. The early work of Duchon~\cite{duchon1977splines} is a canonical example. Historically, roughness minimization has been achieved through the use of smoothing splines; a thorough exposition of smoothing splines can be found in the book by Gu~\cite{gu2013smoothing}. The application of smoothing splines requires a trade-off between accuracy and roughness minimization, since aggressively penalizing roughness tends to degrade accuracy. Therefore, much work has been devoted to parametrizing this trade-off appropriately. Craven and Wahba~\cite{craven1978smoothing} developed the influential ``cross-validation'' approach, which  is expanded upon by Gu~\cite{gu1992cross}.

The functional approximation used in this article is based on global tensor-product B-splines, but a number of other spline-based regression methods have been proposed. Truncated thin-plate splines were used by Wood~\cite{wood2003thin} to improve the efficiency of thin-plate regression splines while maintaining their characteristic stability. Lee et al.~\cite{lee1997scattered} utilized hierarchies of B-splines to fit unstructured data points, but the instability arising from sparse point distributions was not treated explicitly. Francis et al.~\cite{francis2018scattered} consider a two-step process for resampling unstructured point clouds with variable point density onto unstructured grids. While this method does not construct a functional approximation, it does show good performance as a resampling methodology.

Our novel adaptive regularization procedure was first explored in the dissertation of the second author~\cite{yehDissertation}. The method is also directly inspired by the work of El-Rushaidat et al.~\cite{elrushaidat2021}, in which a two-level regularization process was introduced in the context of resampling unstructured data onto structured meshes. However, their method requires an ad-hoc selection of the criteria to switch between high and low regularization strengths, as well as an application-dependent overall level of smoothing. A notable contribution in our work is a continuously varying regularization strength (not two-level) that is adapted automatically.

\section{Background on B-Splines}\label{sec:background}
In this section, we provide a brief overview of the basic definitions and constructions necessary to describe B-spline models for scientific data. A thorough presentation on the fundamental theory of B-splines can be found in the books by de Boor~\cite{deBoor1978guide} and by Piegl and Tiller~\cite{piegl1997nurbs}.
\subsection{B-Spline Curves}
A one-dimensional B-spline curve of degree $p$ in $\mathbb{R}^D$ is a parameterized curve
\begin{equation}\label{eq:spline-def}
C(u) = \sum_{j=0}^{n-1} N_{j,p}(u) P_j,
\end{equation}
where each $N_{j,p}$ is a piecewise-polynomial function of degree $p$, and each $P_j \in \mathbb{R}^D$ is a ``control point'' in $D$-dimensional space. 

The B-spline basis functions, denoted $N_{j,p}$, are defined on the parameter space $[0,1] \subset \mathbb{R}$, which is divided by a nondecreasing sequence of ``knots'' $t_0 \leq t_1 \leq \ldots \leq t_{n+p} \in [0,1]$. Each basis function $N_{j,p}$ is a bump function in $[t_j, t_{j+p+1}]$ and zero elsewhere.\footnote{%
  We consider only ``clamped'' knot sequences in this paper; thus, the first $p+1$ knots are always $0$ and the last $p+1$ knots are always $1$.%
}
In this paper, we will assume that the degree of the B-spline is fixed and drop the $p$ subscript, instead denoting the $j^{th}$ function as $N_j$.

In order to simplify notation when describing high-dimensional tensor product splines, we use multi-indices to index quantities in multiple dimensions simultaneously.  A multi-index $\alpha = (\alpha^1, \ldots, \alpha^d)$ is a $d$-tuple of nonnegative indices, where the sum of components of $\alpha$ is denoted $\abs{\alpha} = \sum_k \alpha_k$. 

We will often consider index sets for our multi-indices in the form of 
\begin{equation}\label{eq:index-set}
  A = \{\text{all } \alpha \in \mathbb{N}^d \text{ such that } 0 \leq \alpha^k < n_k \text{ for } 1 \leq k \leq d\},
\end{equation}
where $n_k$ are previously defined positive numbers. We impose a lexicographic ordering on these sets, 
and denote by $\lidx{\alpha}_A$ the index of $\alpha$ in the lexicographic ordering of $A$.
In the following sections, we consider matrices in which each column corresponds to a multi-index. In this scenario, we list multi-indices in lexicographic order; thus, the multi-index $\alpha$ corresponding to the $j^{th}$ column satisfies $\lidx{\alpha}_A = j$.

\subsection{Tensor Product B-Splines}
Tensor product B-splines are a natural extension of B-spline curves to higher-dimensional manifolds, such as surfaces, volumes, and hypervolumes. Here, we denote by $d$ the dimension of the tensor product volume and $D$ the dimension of the ambient space (for instance, a 2D surface in 3D space would correspond to $d=2$, $D=3$). The parameter space for a $d$-dimensional tensor product B-spline is $[0,1]^d$, which is divided by $d$ different knot vectors $\mathbf{t}_k = \{t_k^j\}_{j=0}^{n_k + p}$, $k=1,\ldots,d$. 

Given a tuple $u = (u^1, \ldots, u^d) \in [0,1]^d$, the tensor product basis functions are defined as
\begin{equation}
  N_\alpha(u) = \prod_{k=1}^d N^k_{\alpha^k}(u^k).
\end{equation}
where $\alpha$ is a multi-index as described above and $N_{\alpha^k}^k$ is the $(\alpha^k)^{th}$ basis function with respect to the knot vector $\mathbf{t}_k$.
With $n_k + p + 1$ total knots in each dimension, there are $n_k$ basis functions in each dimension.\footnote{%
Here we assume for simplicity that the degree of the B-spline is the same in each dimension, but the degree can vary in practice if desired.}
Therefore, the total number of tensor product basis functions is $n_{tot} = \prod_{k=1}^d n_k$, which is also the total number of control points for the the tensor product spline.

A $d$-dimensional tensor product B-spline in $\mathbb{R}^D$ is a function of the form
\begin{equation}\label{eq:tensor-spline}
  C(u) = \sum_{\alpha \in A} N_\alpha(u) P_\alpha,
\end{equation}
where $A$ is the set of all basis functions and $P_\alpha \in \mathbb{R}^D$ for each $\alpha$.

\subsection{Optimal Control Points}
Given a collection of knot vectors and polynomial degree, the best-fit B-spline to a given data set is determined by a linear least-squares minimization problem. Let $\{Q_i\}_{i=0}^{m-1}$ be the list of points in $\mathbb{R}^D$ to be approximated with a $d$-dimensional tensor product spline. For each $0 \leq i < m$, let $v_i \in [0,1]^d$ be the parameter tuple corresponding to the point $Q_i$. The optimal control points are determined by the least-squares minimization problem:
\begin{equation}\label{eq:standard-min}
  \{\hat{P}_j\} = \argmin_{P_j} \sum_{i=0}^{m-1} \norm{Q_i - C(v_i)}^2.
\end{equation}

This minimization problem can be rewritten in normal form by differentiating the objective function in Equation~\eqref{eq:standard-min} with respect to each of the control points. The normal system reduces to the matrix equation
\begin{equation}
\mathbf{N}^T \mathbf{N} \mathbf{P} = \mathbf{N}^T \mathbf{Q},
\end{equation} where
\begin{equation}\label{eq:basic-mat-defs}
  \begin{aligned}
  &\mathbf{N}_{ij} = N_\alpha(v_i) &\text{where } \lidx{\alpha}_A = j\\
  &\mathbf{P}_{ij} = P_\alpha^j, &\text{where } \lidx{\alpha}_A = i\\
  &\mathbf{Q}_{ij} = Q_i^j &
  \end{aligned}
\end{equation}
The superscripts in the above equations index the components of the vectors $P_\alpha$ and $Q_i$.
$\mathbf{N}$ is a $m \times n_{tot}$ matrix, often called the ``B-spline collocation matrix,'' $\mathbf{P}$ is an $n_{tot} \times D$ matrix with each row containing a control point, and $\mathbf{Q}$ is a $m \times D$ matrix with each row containing an input point.

Typically, the matrix $\mathbf{N}$ is very sparse, and this system may be solved with an iterative method or sparse direct solver. However, as we show in the following section, this system is ill-conditioned when the sample density of the input points $P_i$ varies from region to region.

\section{Adaptive Regularization}\label{sec:method}
\begin{figure*}[t]
  \centering
  \includegraphics[width=165pt, height=90pt]{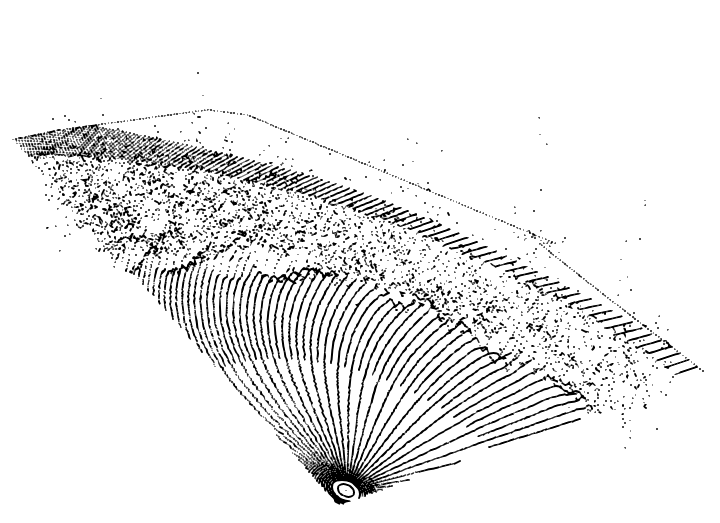}
  \includegraphics[width=165pt, height=90pt]{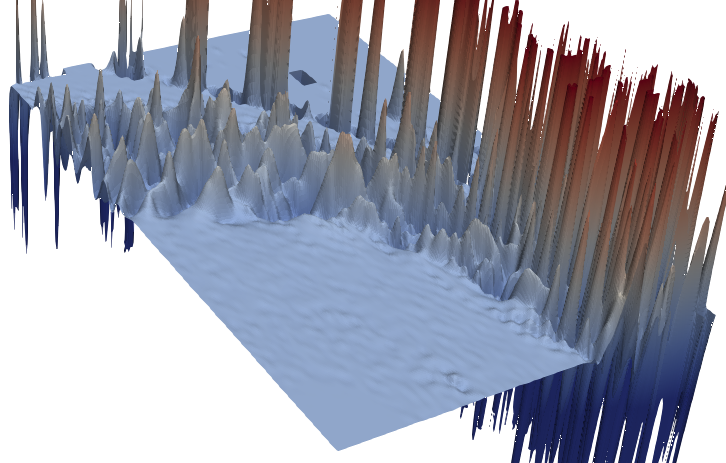}
  \includegraphics[width=165pt, height=90pt]{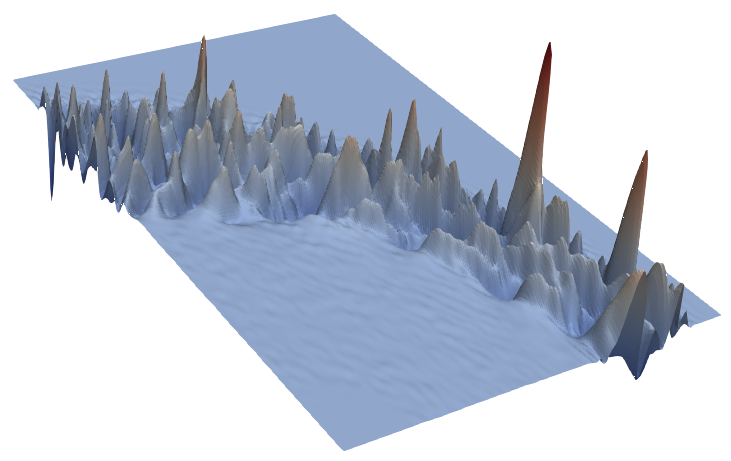}
  \caption{Left: A data set with nonuniform point density. Center: Best-fit B-spline model without regularization. Right: B-spline model with adaptive regularization. The center image is cropped; spike artifacts in this model extend well outside the frame.}
  \label{fig:xgc-sample}
\end{figure*}

A significant challenge when modeling unstructured data with tensor product B-splines is the (ill-)conditioning of the fitting procedure. Generally speaking, tensor product B-spline models can oscillate strongly due to overfitting in regions where input data is sparse (see Figure~\ref{fig:xgc-sample}).
Our method of adaptive regularization produces a unique solution to systems which would otherwise be rank-deficient and improves the overall conditioning of the system. In practice, the adaptively regularized models possess fewer oscillatory artifacts and do not exhibit any divergent behavior in our testing. In contrast to standard regularization techniques, our method does not smooth out the model indiscriminately -- instead, it regularizes only those regions that require smoothing.

This technique employs a spatially-varying regularization strength that is computed automatically as a function of the relative positioning of input data points to the B-spline knots. In general, the regularization strength increases in regions with little to no input data and decreases (potentially to zero) in regions ``saturated'' with input points. When the regularization strength is zero throughout a region of the domain, no smoothing is performed in that region; therefore, any sharp features present in densely sampled regions of the domain will be preserved by the adaptive regularization procedure.

\subsection{Second Derivative Regularization}\label{ssec:method-second-derivs}
Standard roughness minimization can be formulated as a penalized least-squares minimization problem, similar to Equation~\eqref{eq:standard-min}. The penalty term weights the size of the second derivative at each point with a new parameter, which we denote as $\lambda > 0$. The control points of the regularized spline are defined by:
\begin{equation}\label{eq:penalized-min}
  \{\hat{P}_j\}  = \argmin_{P_j} \left(\sum_{i=0}^{m-1} \norm{Q_i - C(v_i)}^2 + \lambda^2 S(C) \right),
\end{equation}
where $S(C)$ approximates the size of the second derivatives of $C$.

Let $w_\alpha \in [0,1]^d$ be the parameter that maximizes the value of $N_\alpha$. Let $\partial^\delta$ denote the partial derivative where the order of derivative in each dimension is given by the components of multi-index $\delta$.\footnote{%
For example, $\partial^{(2,0)} f = \partial^2 f / \partial x_1^2$, and $\partial^{(0,2)} f = \partial^2 f / \partial x_2^2$, while $\partial^{(1,1)} f = \partial^2 f / (\partial x_1 \partial x_2)$.
} 
We define $S(C)$ to be 
\begin{equation}\label{eq:deriv-approx}
  S(C) = \sum_{\alpha \in A} \sum_{\abs{\delta}=2} \norm*{\partial^{\delta}C(w_\alpha)}^2.
\end{equation}
Note that the summation above is a sum over all derivatives of order 2, including mixed partial derivatives.

Equation~\eqref{eq:penalized-min} can be converted into a system of equations in the same way as Equation~\eqref{eq:standard-min}. The only additional step is computing the derivative of $S(C)$ with respect to the control points. 
In matrix form, the new system is 
\begin{equation}
  \begin{pmatrix} \mathbf{N}^T & \lambda \mathbf{M}^T \end{pmatrix}
  \begin{pmatrix} \mathbf{N} \\ \lambda \mathbf{M} \end{pmatrix} \mathbf{P} = 
  \mathbf{N}^T \mathbf{Q}
\end{equation}
where
\begin{equation}\label{eq:augment-mat}
  \mathbf{M} = 
  \begin{pmatrix} \mathbf{M}_{\delta_1} \\ \vdots \\ \mathbf{M}_{\delta_k} \end{pmatrix}
  , \quad\text{and } (\mathbf{M}_{\delta})_{i,j} = \partial^\delta N_\beta(w_\alpha),
\end{equation}
for $\lidx{\alpha}_A = i$, $\lidx{\beta}_A = j$. Intuitively, each column of $\mathbf{M}_\delta$ describes the $\partial^\delta$ partial derivative of an individual B-spline basis function. The matrix $\mathbf{M}$ is the concatenation of all the individual $\mathbf{M_\delta}$ matrices, where $\abs{\delta} = 2$. $\mathbf{N}$, $\mathbf{P}$, and $\mathbf{Q}$ are defined as in Section~\ref{sec:background}. 

The novel improvement of our adaptive regularization scheme is to modify the above system of equations by varying the size of $\lambda$ for each column of $\mathbf{M}$. Since each column of this matrix corresponds to a B-spline basis function and control point, variation in the size of $\lambda$ provides a mechanism to modify the smoothing conditions imposed on each control point of the spline individually.
Due to the local support property of B-splines, setting $\lambda_j^{(2)} = 0$ for control points in a given region ``disables'' the regularization in that region, while still allowing for smoothing to be applied elsewhere.\footnote{
  The superscript ``$(2)$'' on each $\lambda_j$ indicates that this term is applied in the context of second derivatives.} 
Algebraically, we replace the scalar parameter $\lambda$ by a diagonal matrix $\mathbf{\Lambda} = \diag(\lambda_1^{(2)},\ldots, \lambda_{n_{tot}}^{(2)})$, where each $\lambda_j^{(2)} \geq 0$, and consider the new linear system
\begin{equation}\label{eq:penalized-matrix}
  \begin{pmatrix} \mathbf{N}^T & (\mathbf{M}\mathbf{\Lambda})^T \end{pmatrix}
  \begin{pmatrix} \mathbf{N} \\ \mathbf{M}\mathbf{\Lambda} \end{pmatrix} \mathbf{P} = 
  \mathbf{N}^T \mathbf{Q}.
\end{equation}
The value of each $\lambda_i^{(2)}$ is computed automatically as a function of the relative positioning between input data points and B-spline knots. 

To better control this function, we introduce a user-specified parameter called the ``regularization threshold,'' denoted $s^*$. Changing the regularization threshold adjusts the criterion by which some regions of the domain are smoothed and others are not. As $s^*$ increases, smoothing constraints will be applied to larger and larger regions in the domain.

Let $s_j$ denote the $j^{th}$ column sum of $\mathbf{N}$ and $\widetilde{s}_j$ the $j^{th}$ column sum of $\mathbf{M}$. Given $s^* \geq 0$, we define the \emph{second derivative local regularization strength} to be
\begin{equation}\label{eq:lambda-j}
  \lambda_j^{(2)} = \frac{\max(s^* - s_j,0)}{\widetilde{s}_j}.
\end{equation}
Thus, $\mathbf{\Lambda}$ is defined such that every column sum of $\binom{\mathbf{N}}{\mathbf{M\Lambda}}$ is no less than $s^*$.

Adapting the regularization strengths $\lambda_j$ this way has a number of important results. 
If $s^* = 0$, then $\mathbf{\Lambda} = 0$ and the minimization becomes the usual least-squares problem. When $s^*$ is small, $\lambda_j$ will be zero unless the $j^{th}$ column sum of $\mathbf{N}$ is small, which is indicative of an ill-conditioned system. Here, adaptive regularization smooths out only those control points which are poorly constrained. 

This formulation also explains why the adaptive regularization method preserves sharp, densely sampled features while smoothing out oscillatory artifacts. In regions of the domain that are densely sampled, control points will be constrained by many data points and thus the corresponding column sum in $\mathbf{N}$ will be relatively large. By choosing $s^*$ to be sufficiently small, all control points in this region will have a regularization strength of zero; i.e. $\lambda_j = 0$. Therefore, the best-fit spline in this region will not be artificially smoothed.

\subsection{First Derivative Regularization}
The adaptive regularization framework described above is defined in terms of second derivatives but can be easily extended to other orders of derivative. We find that constraining first and second derivatives simultaneously is particularly helpful when modeling data sets with no point values at all in certain regions.  This typically happens when the data represent an object with an interior hole or irregular boundary. 

\begin{figure}
  \centering
  \begin{minipage}{.7\linewidth}
    \includegraphics[width=\linewidth]{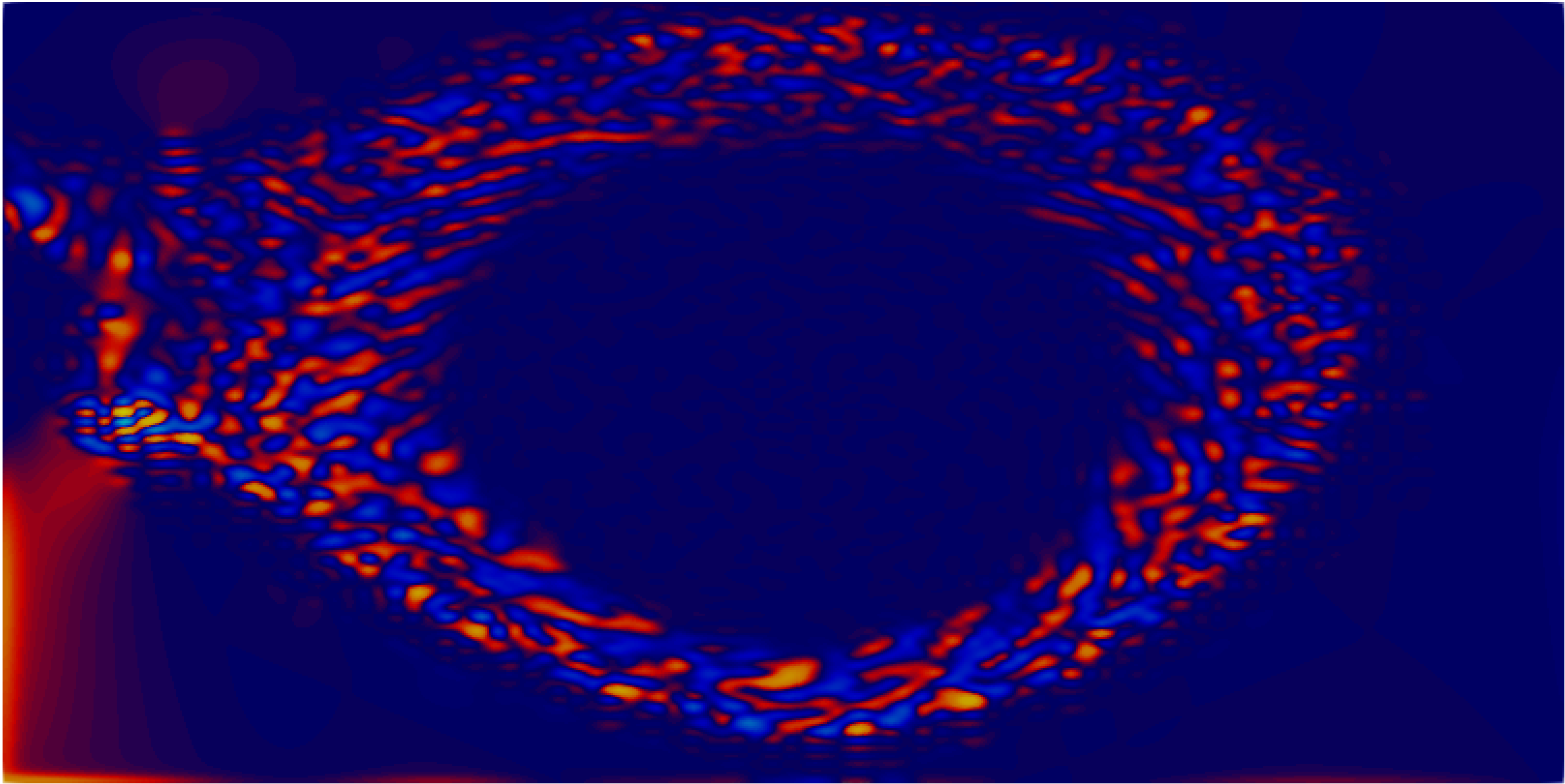}\\
    \includegraphics[width=\linewidth]{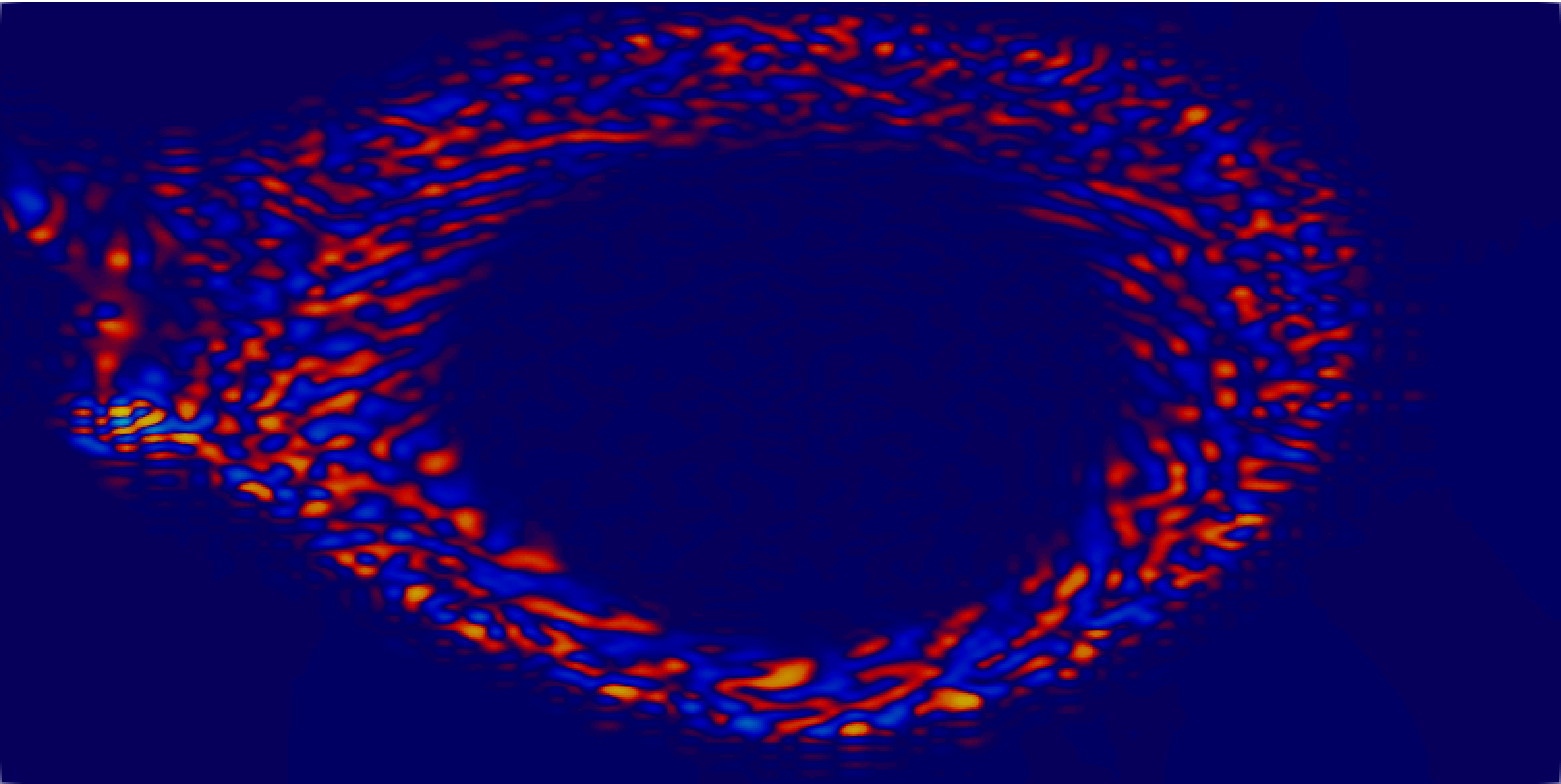} 
  \end{minipage}
  \begin{minipage}{.15\linewidth}
    \includegraphics[width=\linewidth]{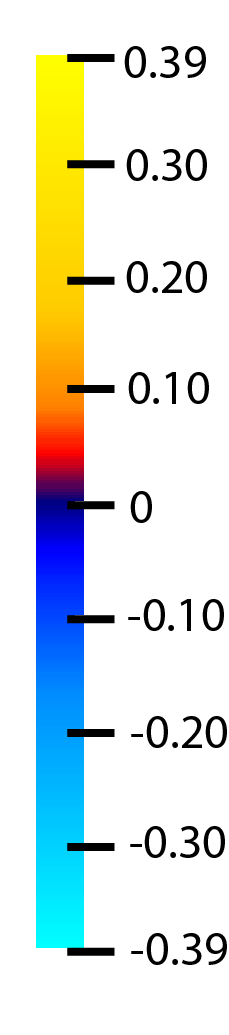}
  \end{minipage}
  \includegraphics[width=.7\linewidth]{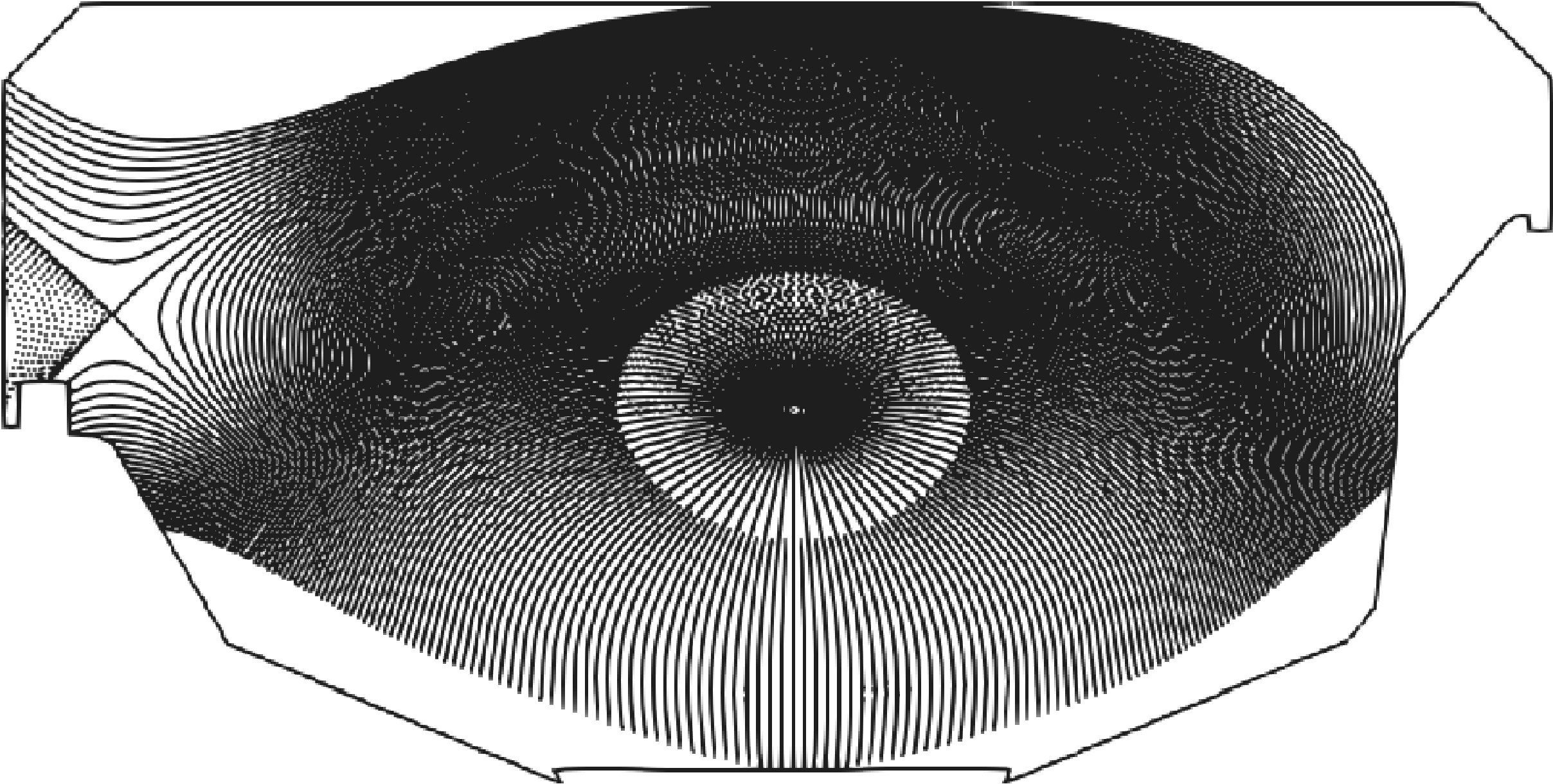}\rule{.15\linewidth}{0pt}\\
  \caption{Effect of using first derivative terms in the regularization scheme. At top: regularization of XGC data set with second derivatives only. At middle: regularization of XGC data set with first and second derivatives. At bottom: point distribution of XGC data set.  See Section~\ref{sec:results} for a description of the data.}
  \label{fig:first-vs-second-derivs}
\end{figure}

In general terms, adaptive regularization with second derivatives produces a model that minimizes curvature in regions where input data is sparse and the fitting problem is ill-conditioned. In regions with no data points, regularized models may still exhibit low-frequency, low-curvature oscillations. As an example of this phenomenon, Figure~\ref{fig:first-vs-second-derivs} shows a B-spline model of a highly non-uniform data set\footnote{
  The data shown is the ``XGC Fusion'' data described in Section~\ref{ssec:datasets}.
}, regularized two ways.  Regularization with second derivatives only (top) produces a distracting ``smearing'' effect on the left edge, where the data set exhibits sharp peaks (note the small yellow regions) in close proximity to a region with no data to constrain the model. In contrast, regularization with first and second derivative terms (middle) can mitigate this numerical artifact. 

In order to include first and second derivative terms in the regularization scheme, the least-squares minimization problem in Equation~\eqref{eq:penalized-matrix} is modified to be
\begin{equation}
  \begin{pmatrix} \mathbf{N}^T & (\mathbf{M_2}\mathbf{\Lambda_2})^T & (\mathbf{M_1}\mathbf{\Lambda_1})^T \end{pmatrix}
  \begin{pmatrix} \mathbf{N} \\ \mathbf{M_2}\mathbf{\Lambda_2} \\
  \mathbf{M_1}\mathbf{\Lambda_1} \end{pmatrix} \mathbf{P} = 
  \mathbf{N}^T \mathbf{Q},
\end{equation}
with
\begin{equation}
\begin{aligned}
  \mathbf{M_2} &= 
  \begin{pmatrix} \mathbf{M}_{\delta_1} \\ \vdots \\ \mathbf{M}_{\delta_k} \end{pmatrix}, \enspace\text{ where } \abs{\delta_i} = 2,\\
  \mathbf{M_1} &= 
  \begin{pmatrix} \mathbf{M}_{\delta'_1} \\ \vdots \\ \mathbf{M}_{\delta'_l} \end{pmatrix}, \enspace\text{ where } \abs{\delta'_i} = 1,
\end{aligned}
\end{equation}
and the matrices $\mathbf{M}_{\delta}$ are defined as in Section~\ref{ssec:method-second-derivs}. Namely, each $\mathbf{N}_{\delta}$ is a matrix of regularization conditions on the $\delta$-partial derivatives. $\mathbf{M_1}$ is the concatenation of all $\mathbf{M}_\delta$ where $\abs{\delta}=1$; in other words, it is the matrix containing regularization conditions on all first-order partial derivatives. Likewise, $\mathbf{M_2}$ is the matrix of regularization terms on all possible second-order partial derivatives.\footnote{$\mathbf{M_2}$ is identical to the matrix $\mathbf{M}$ presented in Section~\ref{ssec:method-second-derivs}, it has simply been relabeled here to match the notation of $\mathbf{M_1}$.}

While first derivative terms are useful for mitigating unwanted visual artifacts, it is typically not appropriate to impose these constraints throughout the domain. Because first derivative regularization penalizes nonzero gradients, such a condition is not consistent with the data when a signal is nonconstant. It is not at all unlikely for a data set to have low point density but model a function with a steep gradient. In such a scenario, the fitting algorithm cannot simultaneously minimize the gradient and accurately model the signal without overfitting the data. Thus, it is important to be selective as to where within a domain first derivative regularization is applied.

Our solution to this problem requires only minor modifications to the method described thus far. Because adding first derivative terms is useful primarily for minimizing oscillations in spatial regions with no input data, we modify the method so that these terms are \emph{active only in such regions}. To achieve this, we define the local regularization strength for first derivative terms to be zero in regions where there  is at least one data point. In analogy to Equation~\eqref{eq:lambda-j}, we define the \emph{first derivative local regularization strength} to be
\begin{equation}\label{eq:first-deriv-reg-strength}
  \mathbf{\Lambda_1} = \begin{pmatrix}\lambda_1^{(1)} & & \\
    & \ddots & \\
    & & \lambda_{n_{tot}}^{(1)}
  \end{pmatrix}, \quad\text{where } \enspace
  \lambda_j^{(1)} = 
  \begin{dcases}
    0 & \text{if } s_j > 0\\
    \frac{s^*}{\widetilde{s}_j} & \text{if } s_j = 0
  \end{dcases},
\end{equation}
and, as before, $s_j$ is the $j^{th}$ column sum of $\mathbf{N}$ and $\widetilde{s}_j$ is the $j^{th}$ column sum of $\mathbf{M_1}$.

This definition warrants some discussion. Recall that the matrix $\mathbf{N}$ is the B-spline collocation matrix (c.f. Equation~\eqref{eq:basic-mat-defs}), with every column representing a basis function and every row representing an input point. Therefore, the $j^{th}$ column sum of $\mathbf{N}$ is positive if and only if there is at least one input point in the local support of basis function $j$. If $s_j = 0$, then it must be the case that there is no input data to constrain the $j^{th}$ control point. This is precisely the scenario illustrated in Figure~\ref{fig:first-vs-second-derivs} in which regularization with first derivatives is most impactful. Finally, we remark that when $s_j = 0$, the definition in Equation~\eqref{eq:first-deriv-reg-strength} is consistent with the earlier definition of second derivative regularization strength in Equation~\eqref{eq:lambda-j}.

\section{Results}\label{sec:results}
We demonstrate the effectiveness of our method with a series of numerical experiments. First, we compare adaptive regularization to uniform regularization where the smoothing parameter has been chosen manually. Next, we study the reconstruction of an analytical signal from sparse samples with varying levels of sparsity. 
For each sparsity level, we report the error and condition number for and unregularized and adaptively regularized model.
We then test the performance of adaptive regularization on data sets with no data in certain regions. In these problems, we construct a B-spline model that extrapolates into regions with no pointwise constraints, and check that the adaptive regularization method produces a reasonable result.

\subsection{Data Sets}\label{ssec:datasets}
The performance of the adaptive regularization method was studied on a collection of two- and three-dimensional point clouds with different characteristics. Some data sets were sampled from analytical functions so that we could compute pointwise errors relative to a ground truth, while other data sets were generated by scientific experiments and simulations.

\textit{2D Polysinc.}
The polysinc\footnote{%
We consider the unnormalized sinc function: $\sinc(x) = \sin(x)/x$, with $\sinc(0)=1$.} data set is a two-dimensional point cloud sampling the function
\begin{equation*}f(x,y) = \sinc\left(x^2 + y^2\right) \sinc\left(2(x-2)^2 + (y+2)^2\right).
\end{equation*}
360,000 point locations are uniformly sampled from the box domain $[-4\pi,4\pi] \times [-4\pi,4\pi]$, except at four disk-shaped regions where the sample rate is 50$\times$ lower. 

\begin{figure}[h]
  \centering
  \includegraphics[width=.6\linewidth]{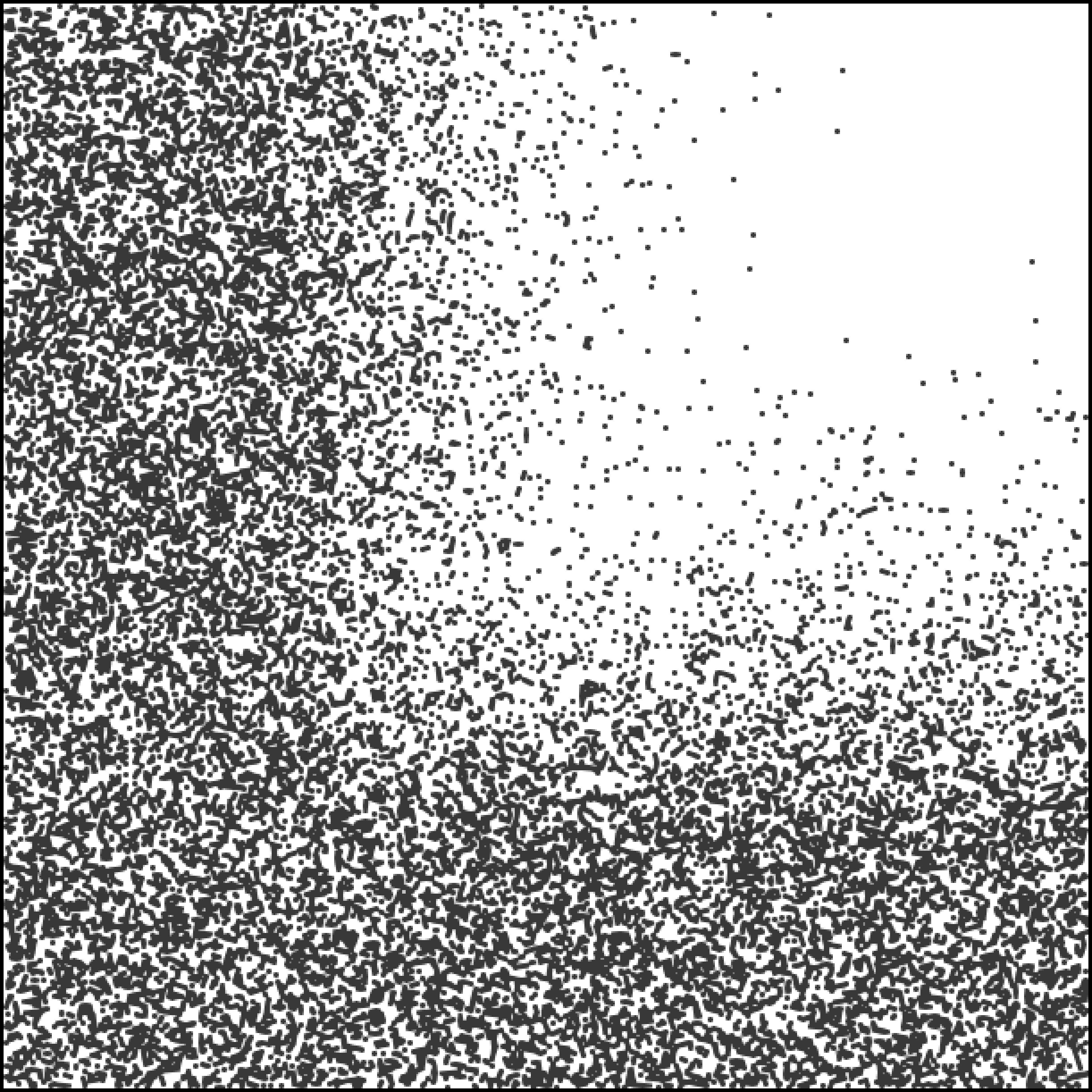}
  \caption{Distribution of points in the modified polysinc data set.}
  \label{fig:polysinc-points}
\end{figure}

\textit{2D Polysinc (Modified).}
The modified polysinc samples the function $f(x,y)$ described above on a non-uniformly distributed point cloud. The data set consists of 22,500 
locations in the box $[-4\pi,4\pi] \times [-4\pi,4\pi]$, where the $(+,+)$ quadrant has the lowest point density and the $(-,-)$ quadrant has the highest point density. Figure~\ref{fig:polysinc-points} shows the sampling locations for this data set.

\textit{XGC Fusion.}
The XGC fusion data set represents a normalized derivative of electrostatic potential in a single poloidal plane of a Tokamak fusion simulation. The data set contains 56,980 points with an irregular boundary and was produced by the XGC code~\cite{xgc} in a simulation of the gyrokinetic equations.

\textit{CMIP6 Climate.}
The CMIP6 climate data set represents ocean surface temperature in a projected box region around Antarctica. The data set contains 585,765 points with a large hole (representing Antarctica) in the center and was produced by a Coupled Model Intercomparison Project (CMIP6)~\cite{cmip6} simulation.

\textit{sahex Nuclear.}
The sahex nuclear data set is derived from a simulation of a single nuclear reactor component, produced with the SHARP toolkit~\cite{yu2016sharp}. The three-dimensional data are bounded by a hexagonal prism and point density is coarser in the $z$ dimension than $x$ and $y$. The data set contains 63,048 points.

\subsection{Comparison of Adaptive versus Uniform Regularization}
\begin{figure*}[t]
  \centering
  \begin{minipage}{.38\textwidth}
    \includegraphics[width=\textwidth]{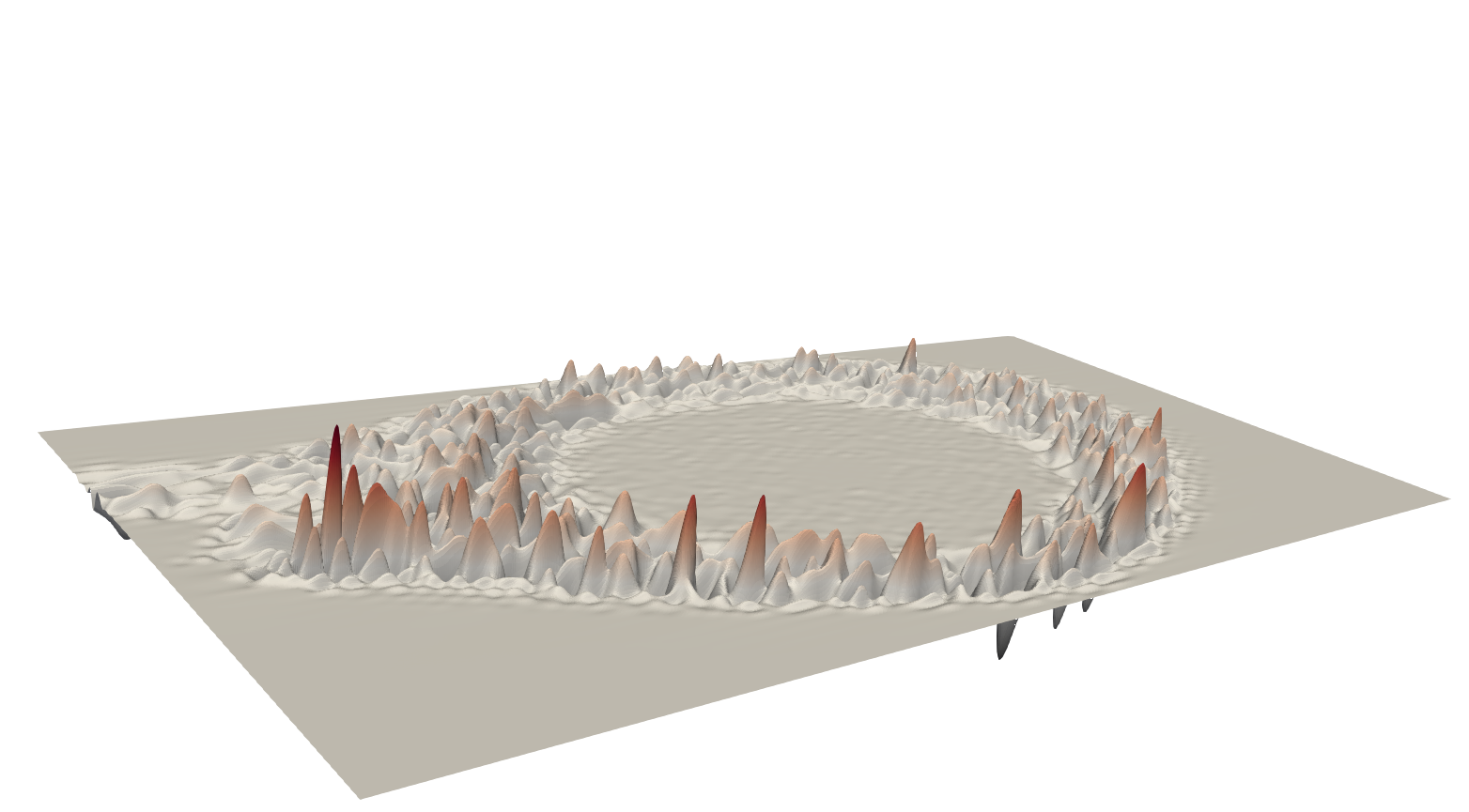}\\[-35pt]
    \includegraphics[width=\textwidth]{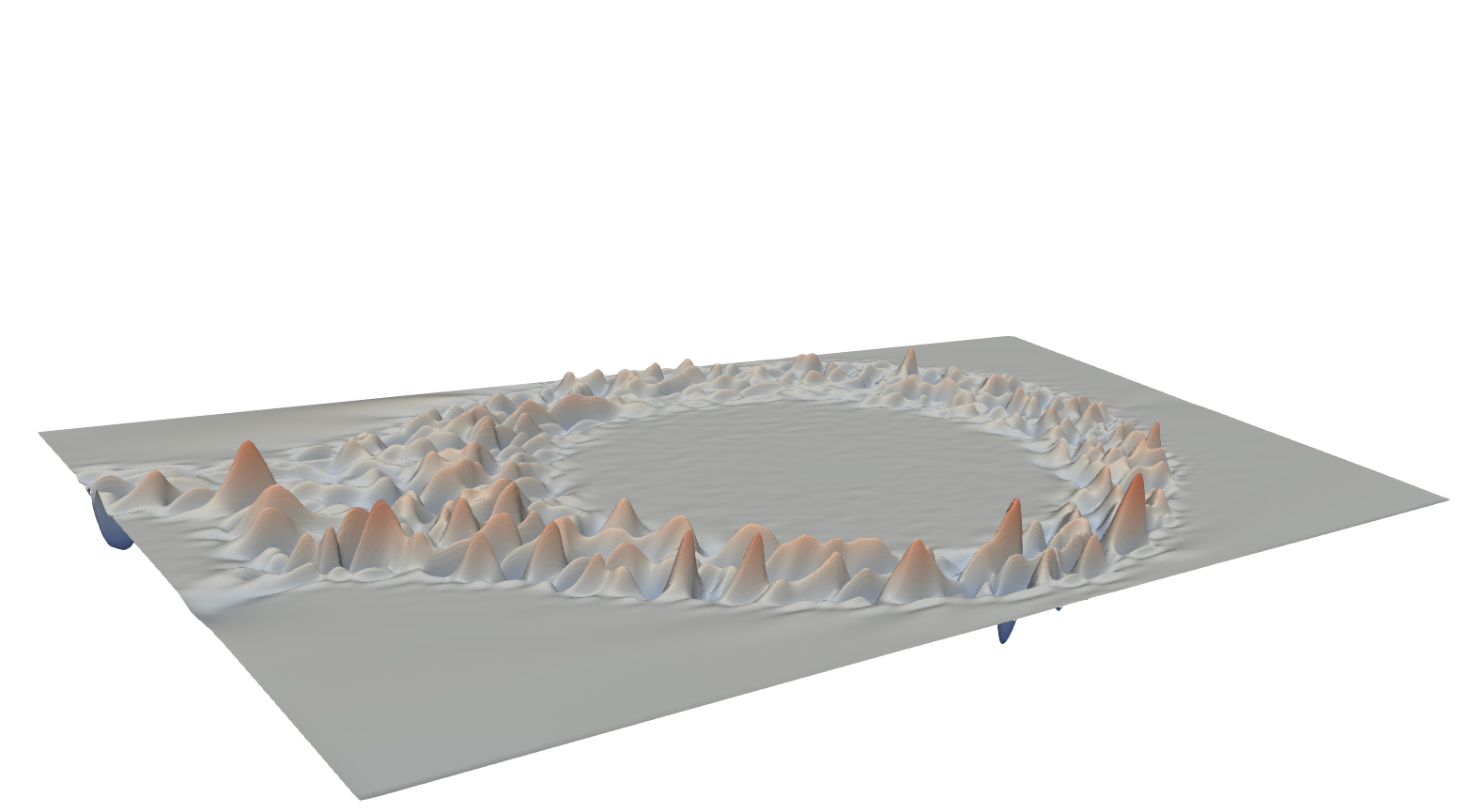}
  \end{minipage}
  \begin{minipage}{.38\textwidth}  
    \includegraphics[width=\textwidth]{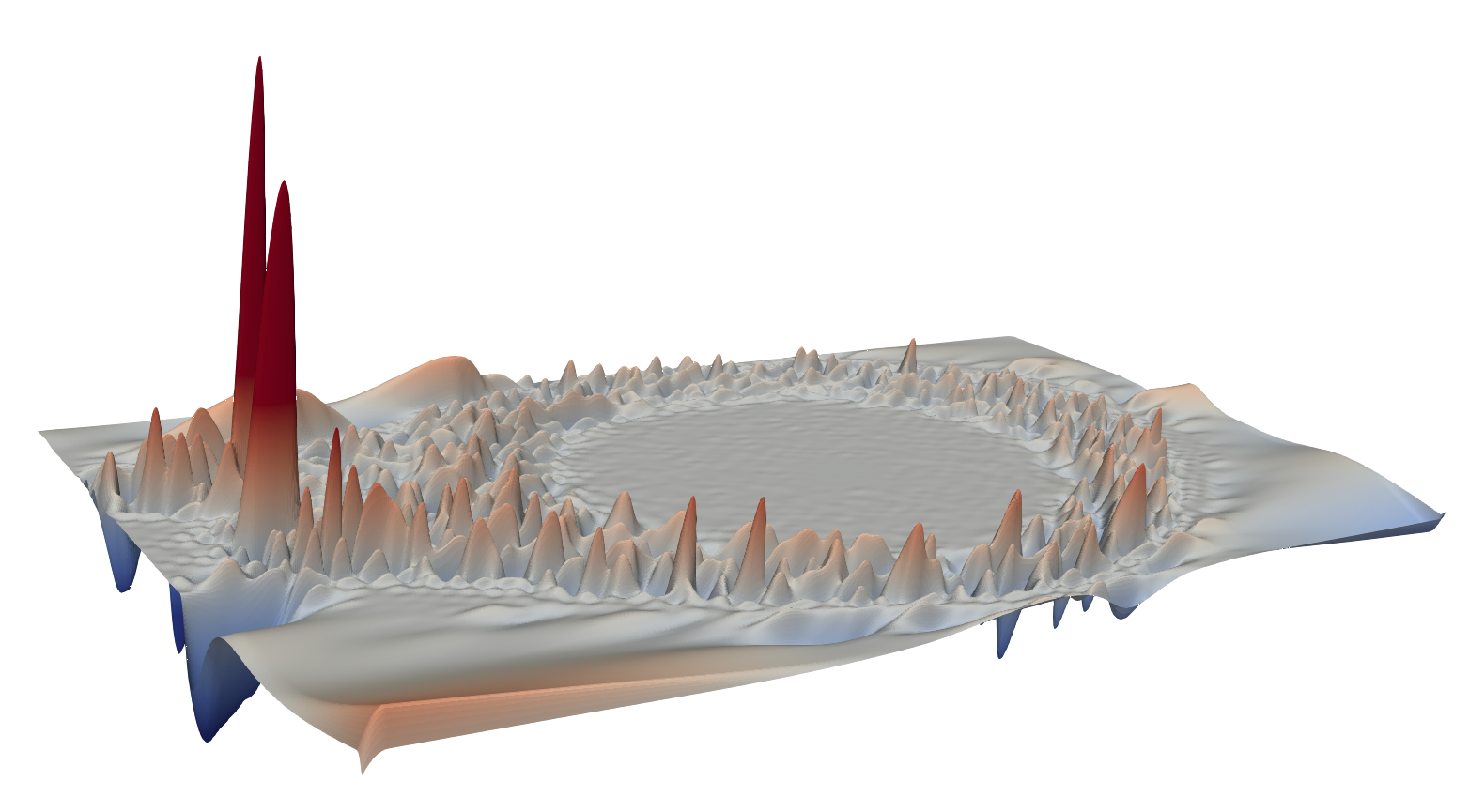}\\[-35pt]
    \includegraphics[width=\textwidth]{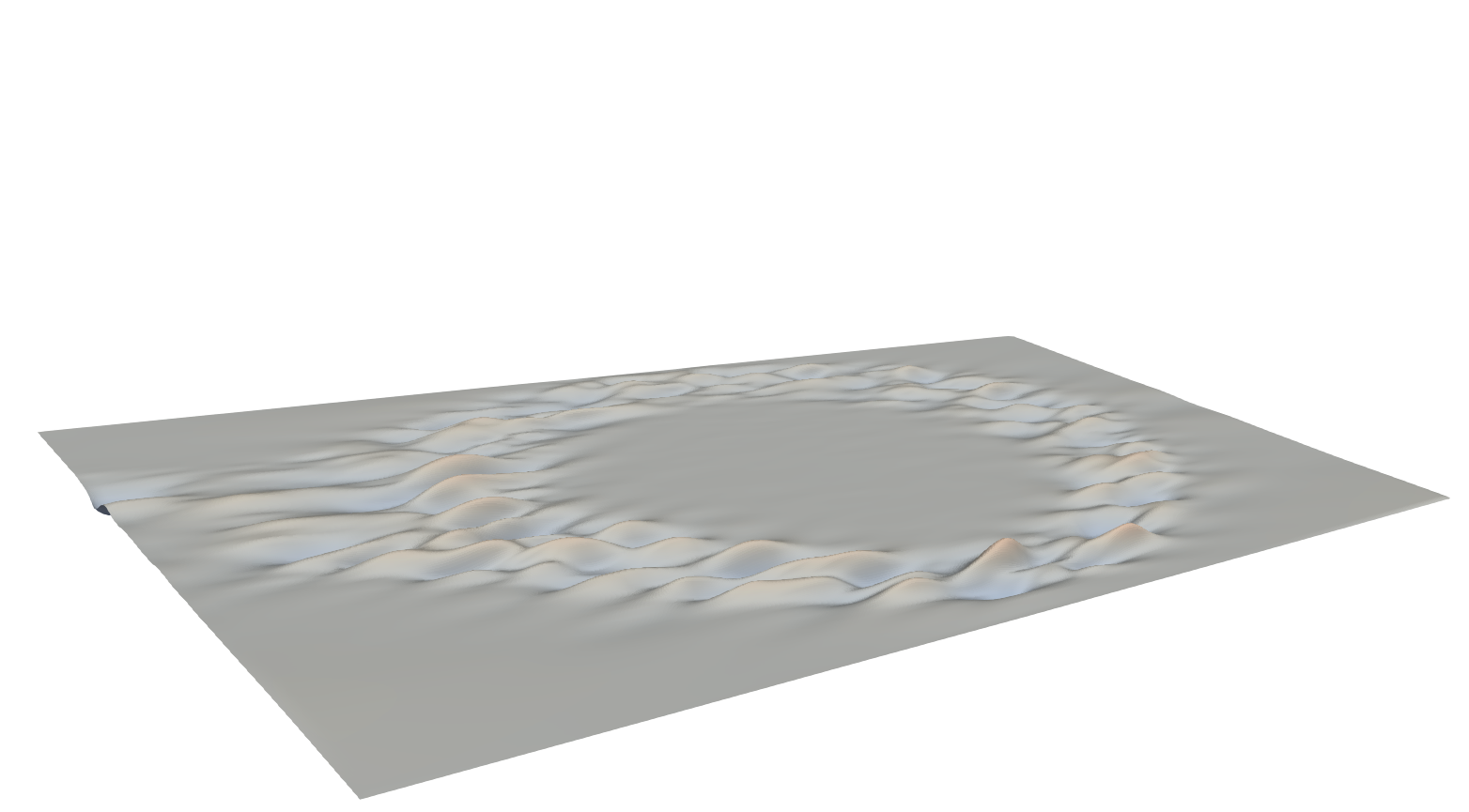}
  \end{minipage}
  \raisebox{-60pt}{\includegraphics[width=30pt]{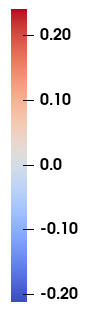}}
  \caption{Comparison of uniform vs adaptive regularization. Clockwise from top-left: Adaptive regularization, uniform regularization with $\lambda = 10^{-6}$, uniform regularization with $\lambda = 10^{-4}$, uniform regularization with $\lambda = 10^{-5}$.}
  \label{fig:uniform-vs-adaptive}
\end{figure*}

Applying a uniform regularization strength to an entire model can produce unsatisfactory results, because sufficiently smoothing oscillatory artifacts can also smooth out sharp features. Figure~\ref{fig:uniform-vs-adaptive} compares our adaptive regularization scheme (with $s^*=6$) against three strengths of uniform regularization. 
The data in Figure~\ref{fig:uniform-vs-adaptive} is the XGC fusion data set, which contains sharp peaks in a ring but is flat inside the ring. Data are sparse or nonexistent outside the ring. The two images at right show a model with uniform regularization that is too weak, causing artifacts (top), or too strong, dampening the features (bottom). The best uniform regularization strength we could find is given at bottom-left, but even in this example the characteristic peaks in the data are smoothed down.

\subsection{Accuracy on Analytical Signals}
\begin{figure*}
  \centering
  \includegraphics[width=.3\textwidth]{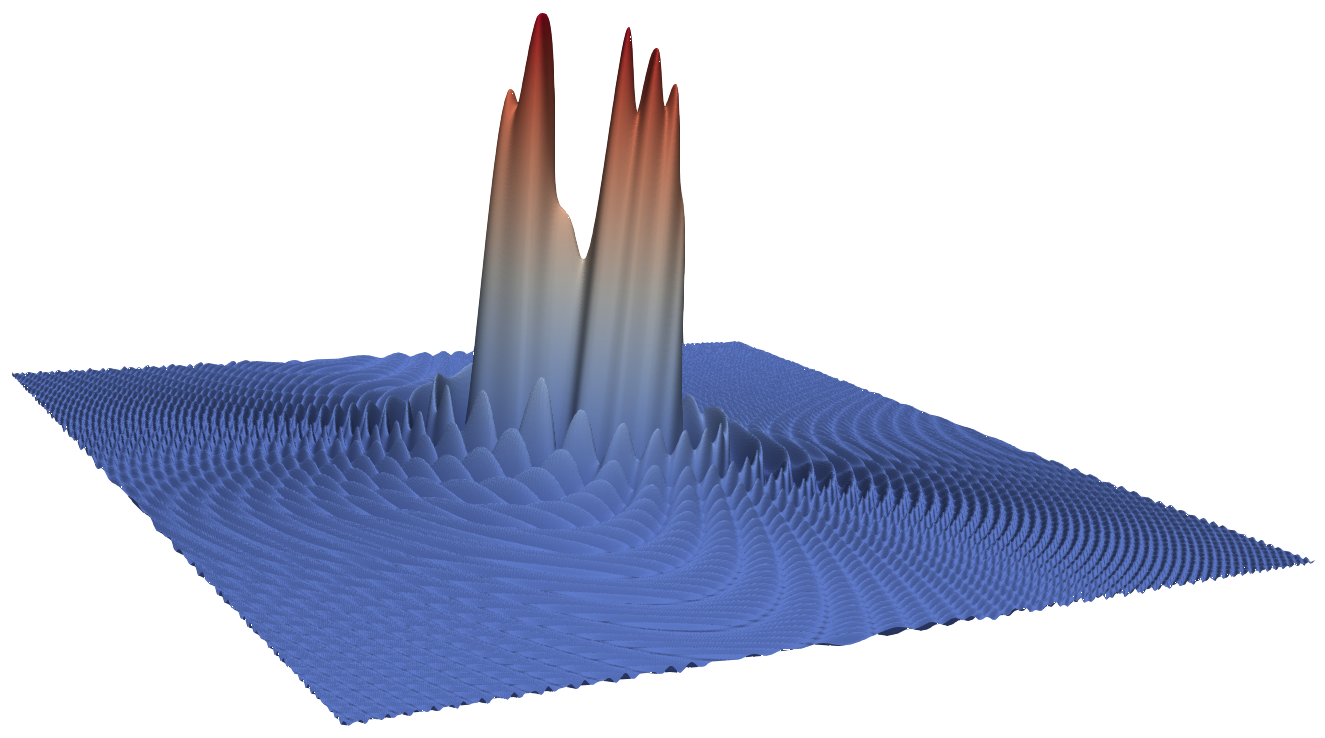} 
  \includegraphics[width=.3\textwidth]{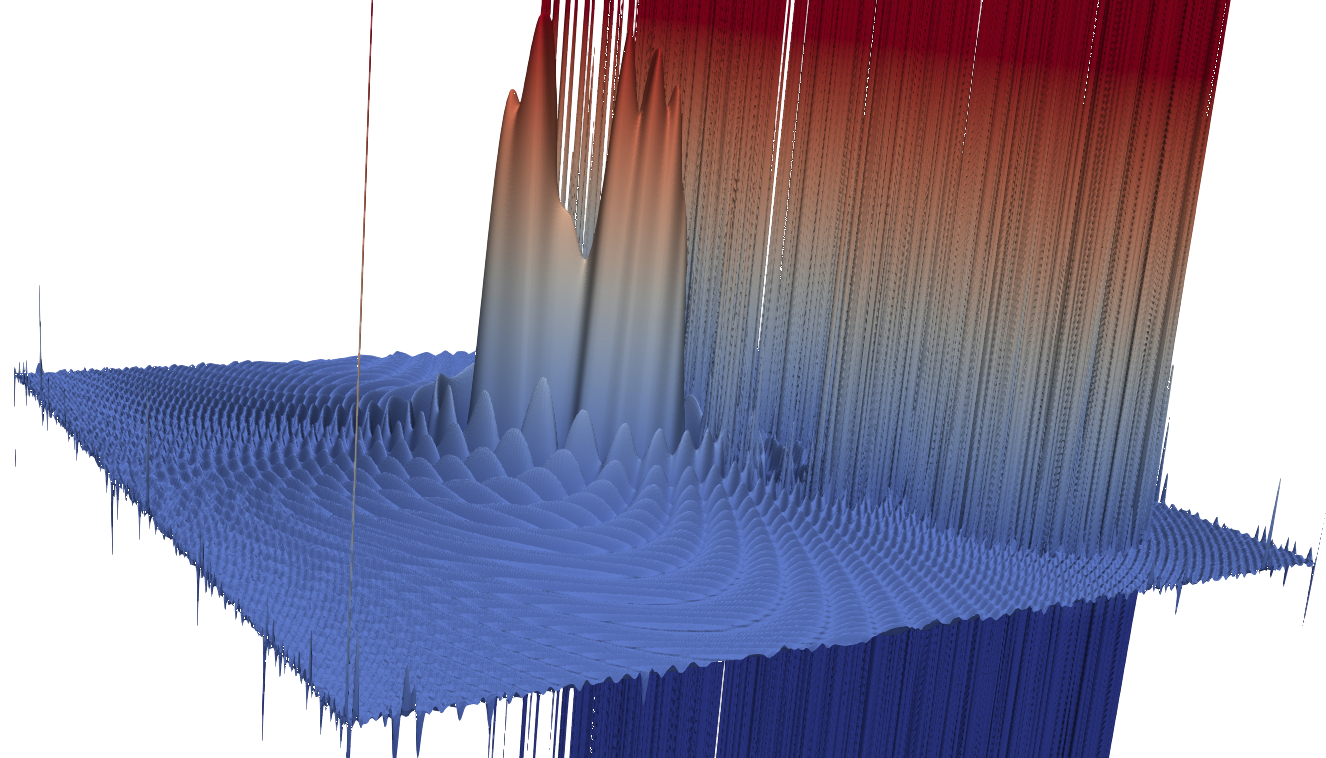}
  \includegraphics[width=.3\textwidth]{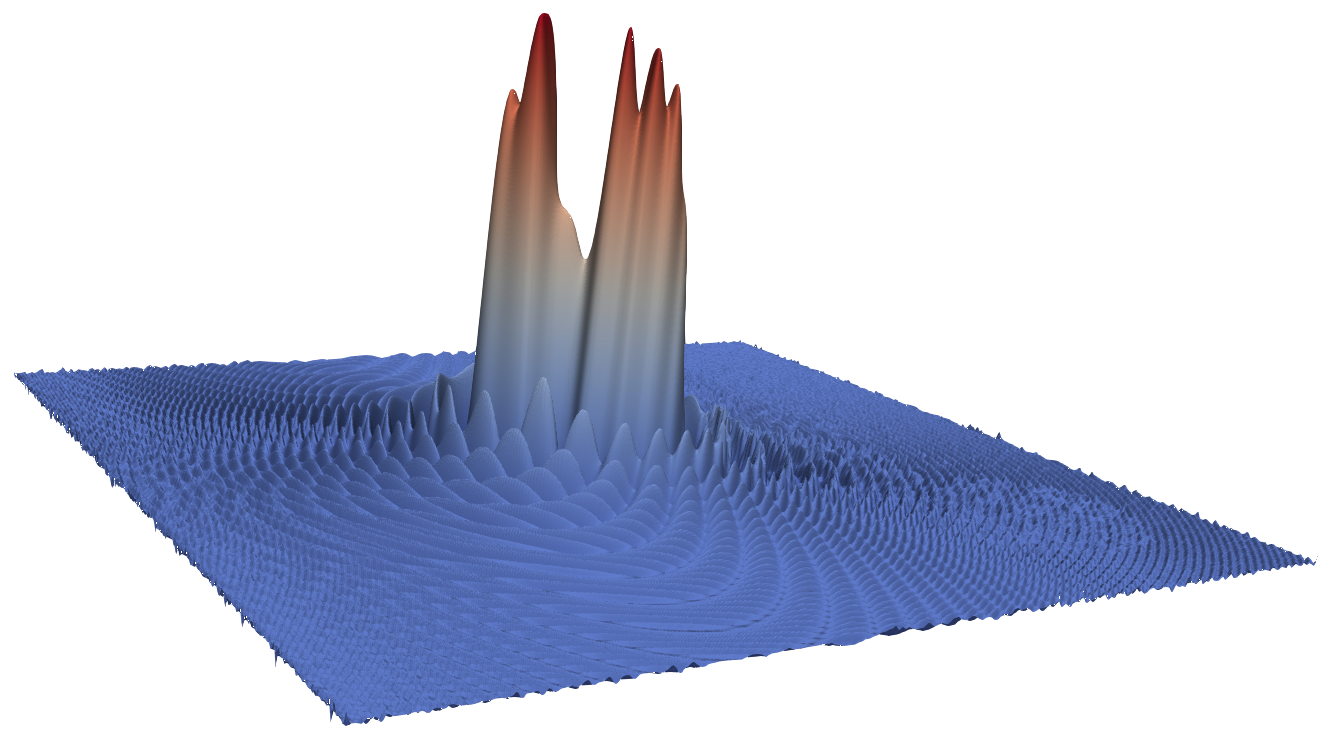}\\
  \hspace{25pt}
  \includegraphics[width=.25\textwidth]{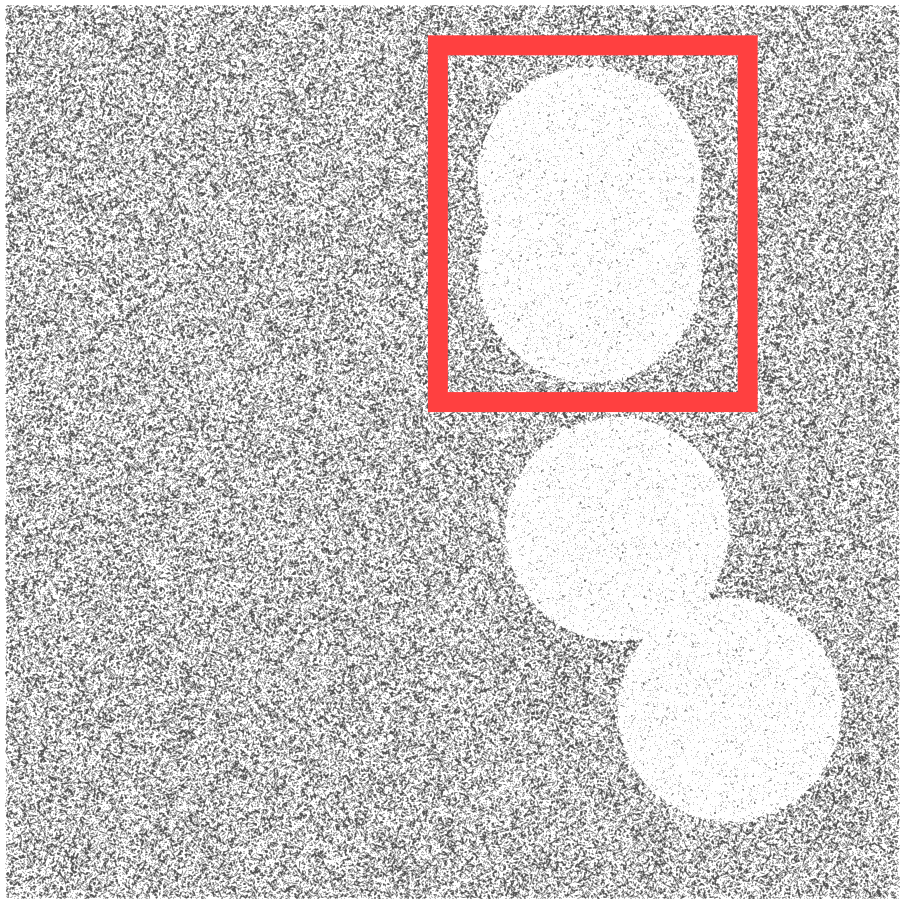}
  \hspace{5pt}
  \includegraphics[width=.25\textwidth]{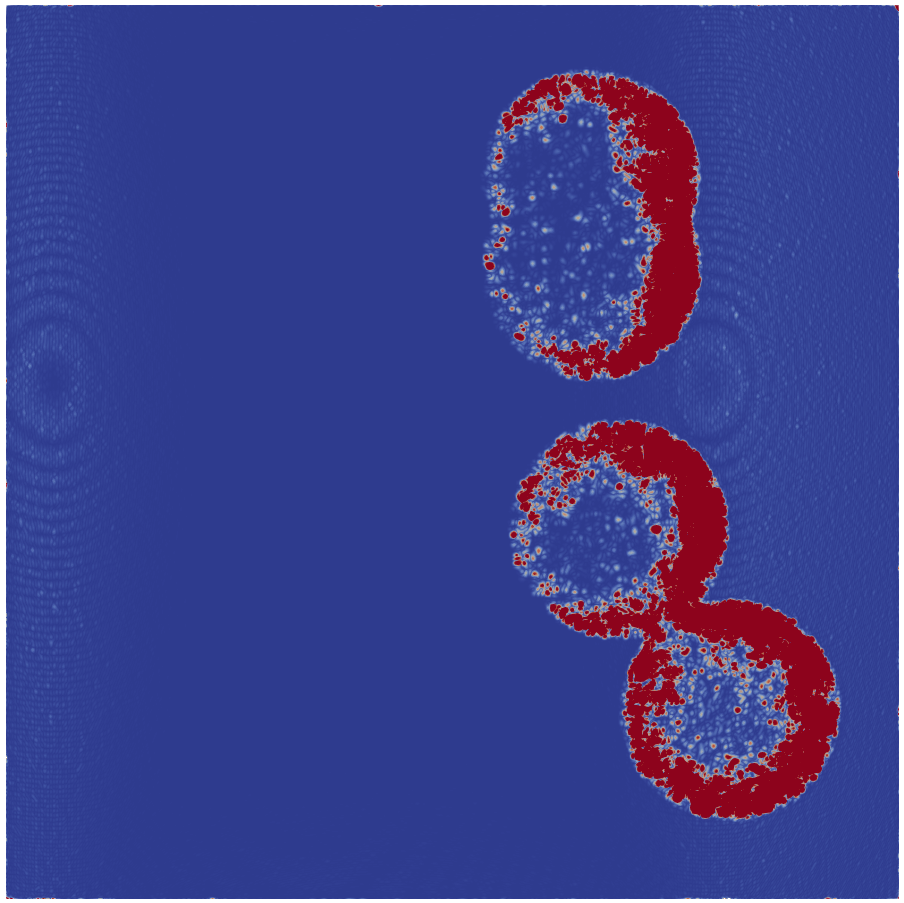}
  \hspace{5pt}
  \includegraphics[width=.25\textwidth]{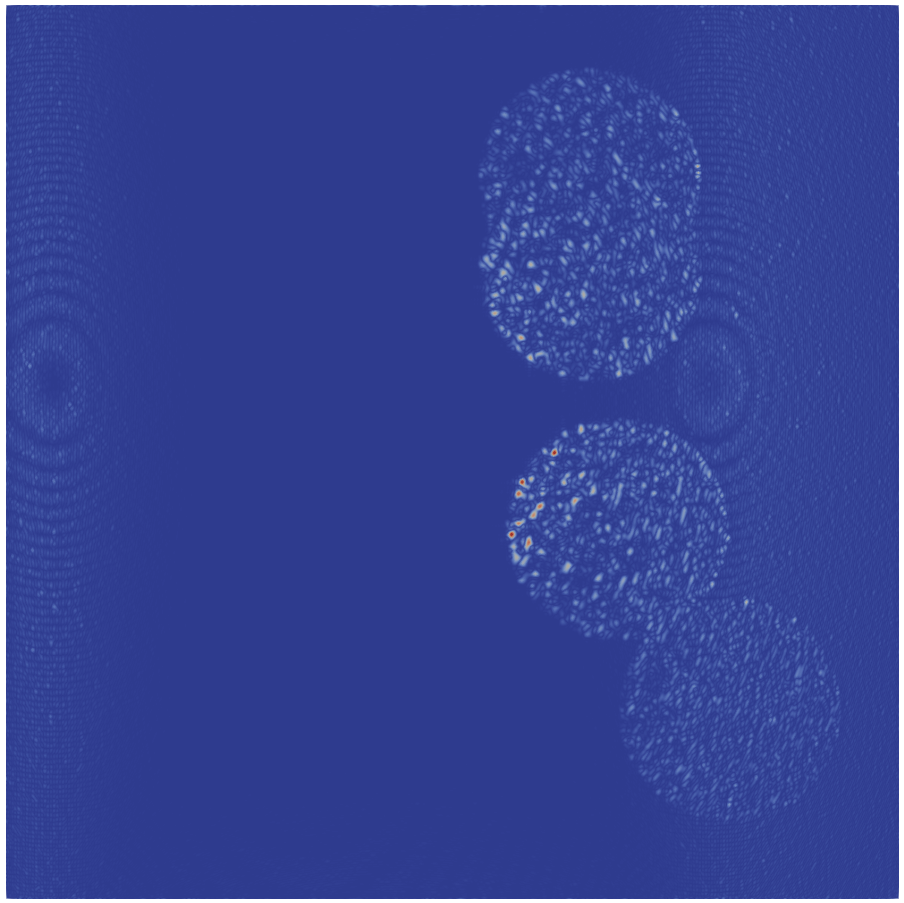}
  \includegraphics[width=.08\textwidth]{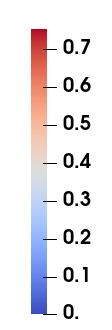}
  \caption{Top row: Synthetic polysinc signal (left), model with no regularization (center), model with adaptive regularization (right). Bottom row: Top down view of input distribution (left), error profile with no regularization (center), error profile with adaptive regularization (right). Area of interest for error calculation is in red at bottom-left.}
  \label{fig:psinc-2d}
\end{figure*}

To quantify the accuracy of B-spline models with adaptive regularization, we consider the oscillatory polysinc function with a highly nonuniform input data set (Figure~\ref{fig:psinc-2d}). We illustrate two B-spline models, one fit without regularization and one with our adaptive regularization ($s^*=1$). Both models are degree four with a 300 $\times$ 300 grid of control points. Without regularization, the model diverges in the regions of low sample density; with adaptive regularization, the model produces an accurate representation even where sample density is low. A top-down view of the error profiles is given in the second row of Figure~\ref{fig:psinc-2d}. A close comparison of the ground-truth (top left) and adaptively-regularized spline (top right) shows that the spline model is not artificially smoothed in dense regions, even though the sparse regions are smoothed.  In particular, our regularization procedure preserves the distinctive oscillations and peaks in the signal.

The degree of sparsity in the input data strongly influences the accuracy of a B-spline model. Table~\ref{tb:error-study} lists the errors in each model for varying levels of sparsity in the voids. The errors are measured in a box around two voids (see Figure~\ref{fig:psinc-2d}) in order to pinpoint the behavior of the models in this region. When the point density is equal inside and outside of the voids (sparsity $=1.0$), error for both models is low. As the voids become more sparse, the error in the unregularized model increases by four orders of magnitude while error in the adaptively regularized model stays essentially flat.

\begin{table*}
  \caption{Model errors as a function of sparseness. Maximum and $L^2$ (average) errors are computed for both adaptively regularized and unregularized models in the vicinity of two voids (see Figure~\ref{fig:psinc-2d}). Condition numbers for both minimization problems are reported at bottom.}
\centering
\begin{tabularx}{\textwidth}{l @{\hspace{5pt}}|@{\hspace{3pt}} 
  *6{>{\centering\arraybackslash}X}}
Sparsity & 0.02 & 0.08 & 0.16 & 0.32 & 0.64 & 1.00\\ \hline
Max Error (reg) & 3.25e-2 & 2.89e-2 & 2.71e-2 & 2.39e-2 & 1.20e-2 & 1.16e-2 \\
Max Error (no reg)  & 1.29e2 & 6.27e2 & 2.36e-1 & 3.65e-2 & 1.20e-2 & 1.16e-2 \\ \hline
$L^2$ Error (reg) & 1.93e-3 & 1.53e-3 & 1.13e-3 & 7.04e-4 & 5.56e-4 & 5.44e-4 \\
$L^2$ Error (no reg) & 2.17e0 & 4.68e0 & 3.42e-3 & 7.35e-4 & 5.56e-4 & 5.44e-4\\ \hline
Condition \# (reg) & 177 & 980 & 289 & 198 & 121 & 189 \\
Condition \# (noreg) & inf & inf & 6.97e4 & 2.58e3 & 1.57e3 & 3.74e3 
\end{tabularx}
\label{tb:error-study}
\end{table*}

Data sparsity also affects the condition number of the least-squares minimization. Table~\ref{tb:error-study} gives the condition number of the matrices $\mathbf{N}$ (`noreg') and $\binom{\mathbf{N}}{\mathbf{M\Lambda}}$ (`reg') for each sparsity level. As sparsity is increased, the condition number of $\binom{\mathbf{N}}{\mathbf{M\Lambda}}$ remains lower and steady but the condition number of $\mathbf{N}$ starts higher and eventually becomes infinite. Condition numbers were computed with the Matlab routine \texttt{svds}. 

\begin{figure*}
\centering
\includegraphics[width=.3\textwidth]{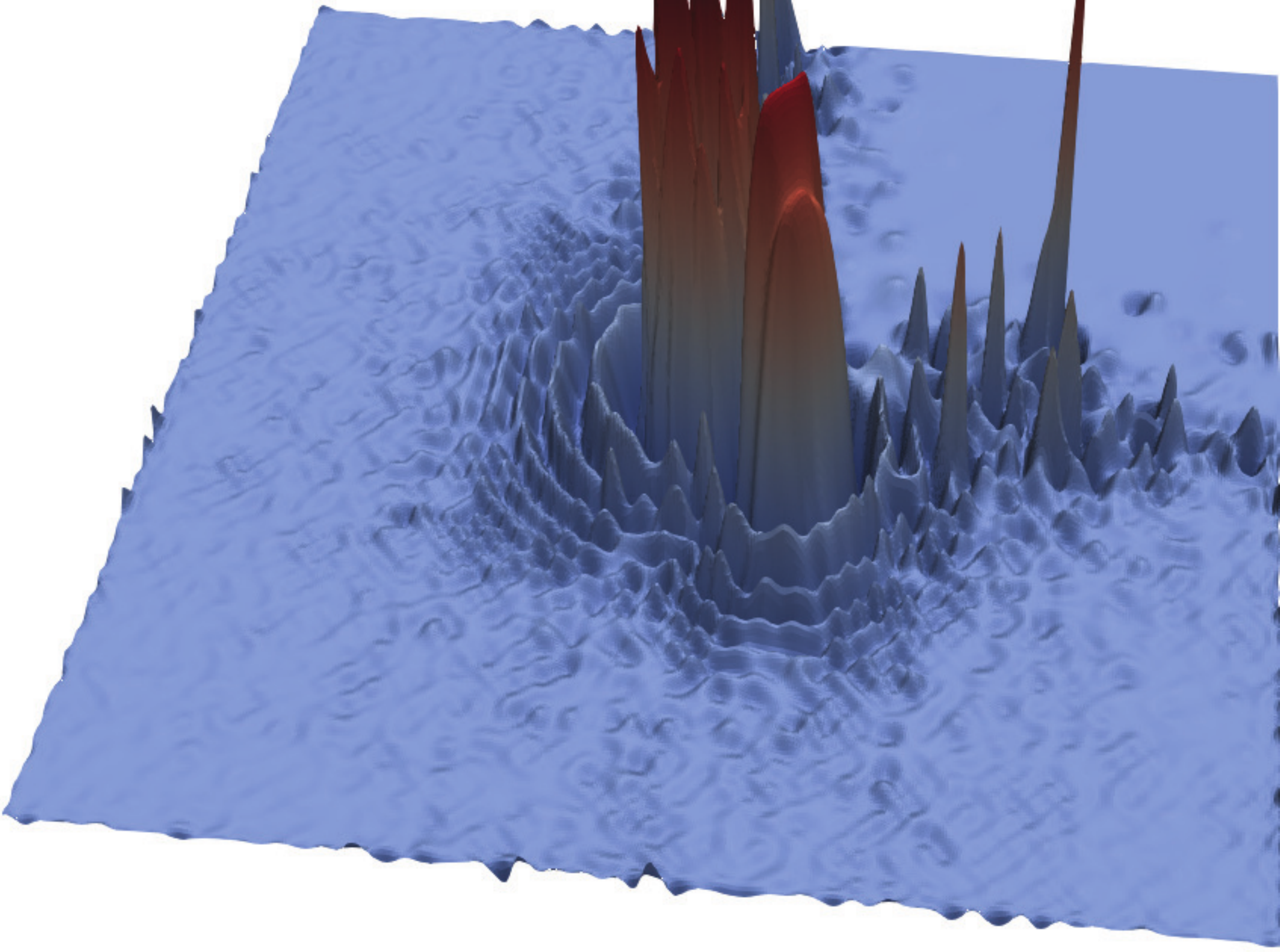}
\includegraphics[width=.3\textwidth]{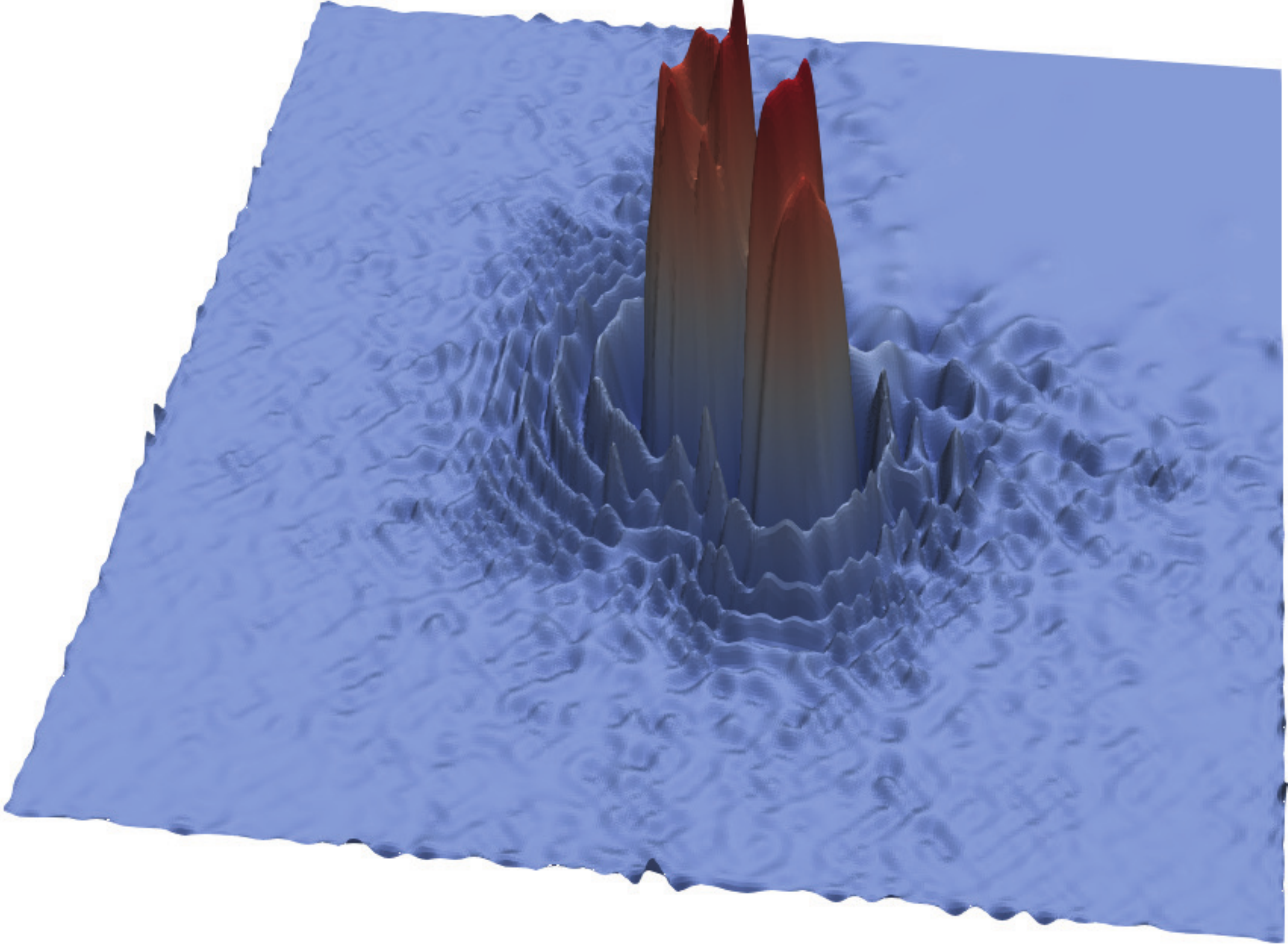}
\includegraphics[width=.3\textwidth]{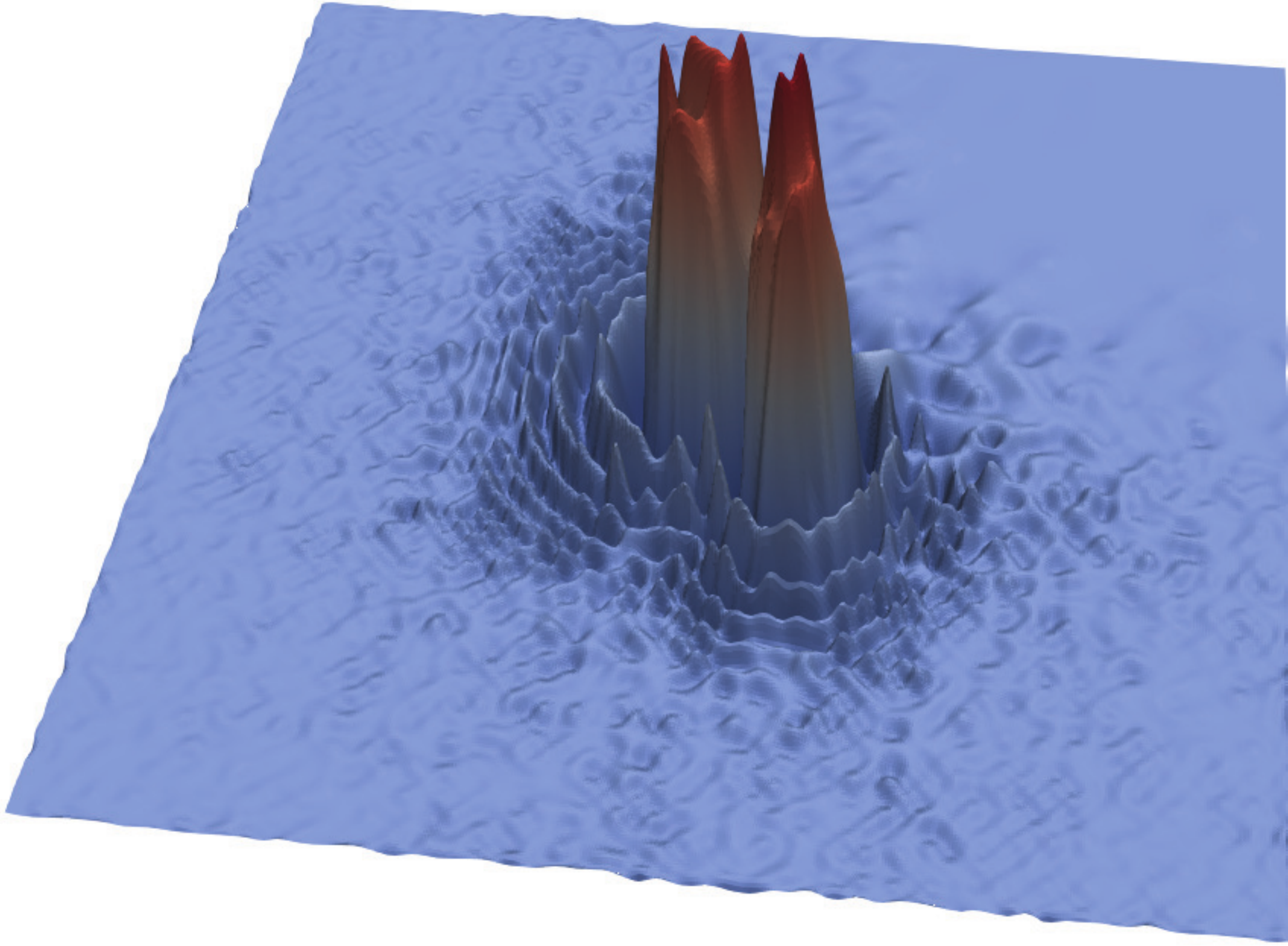}\\
\includegraphics[width=.3\textwidth]{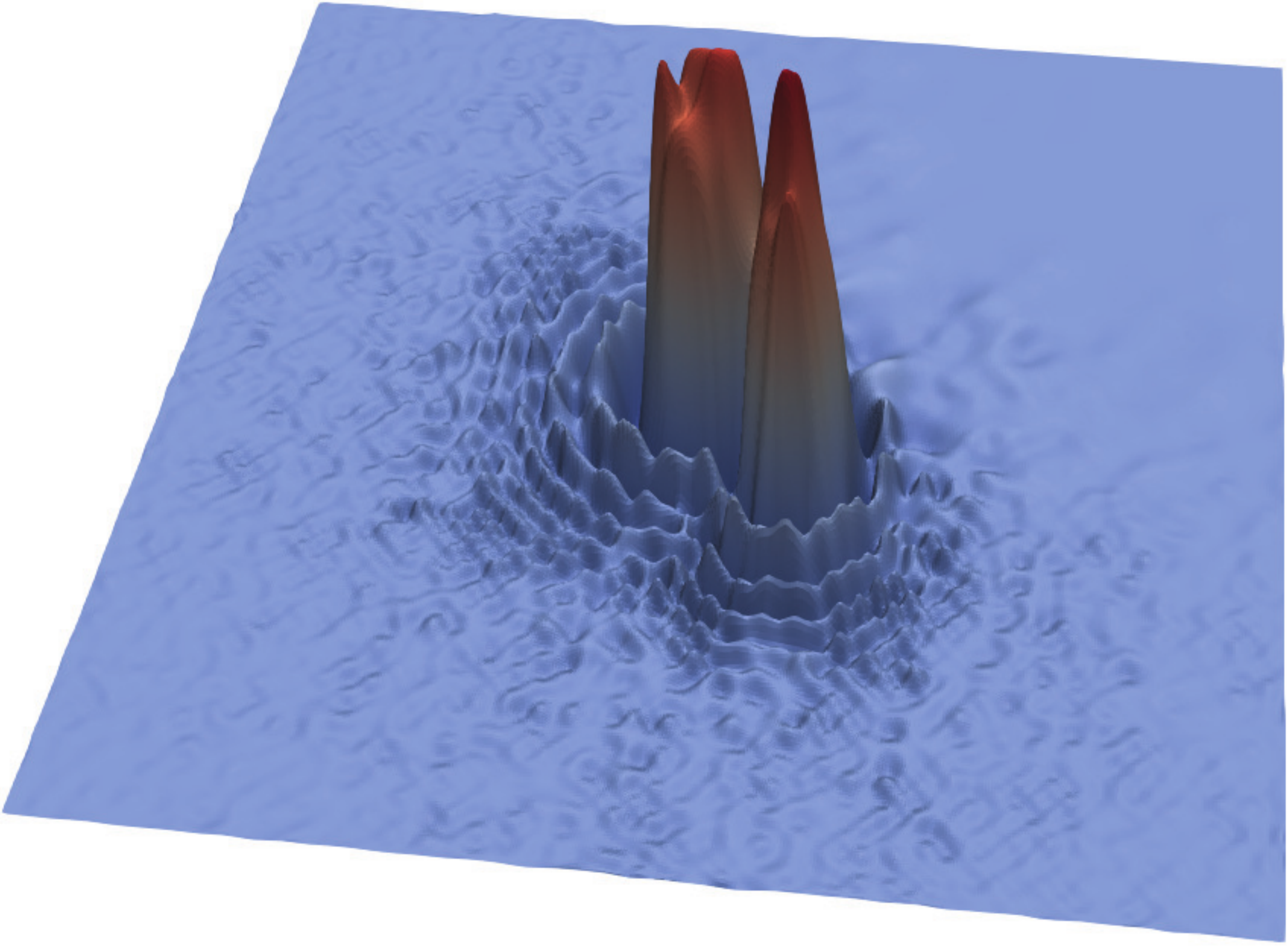}
\includegraphics[width=.3\textwidth]{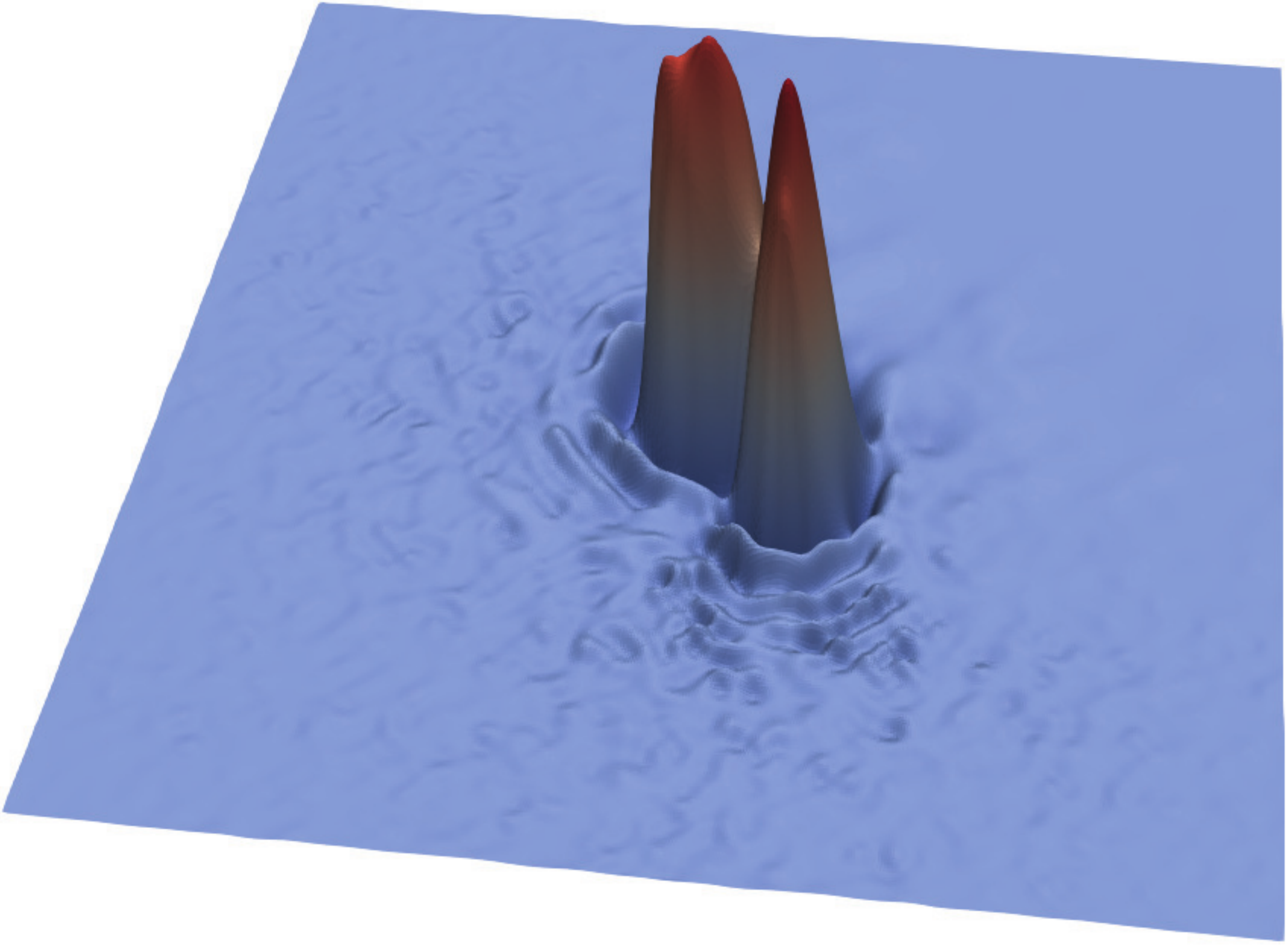}
\includegraphics[width=.3\textwidth]{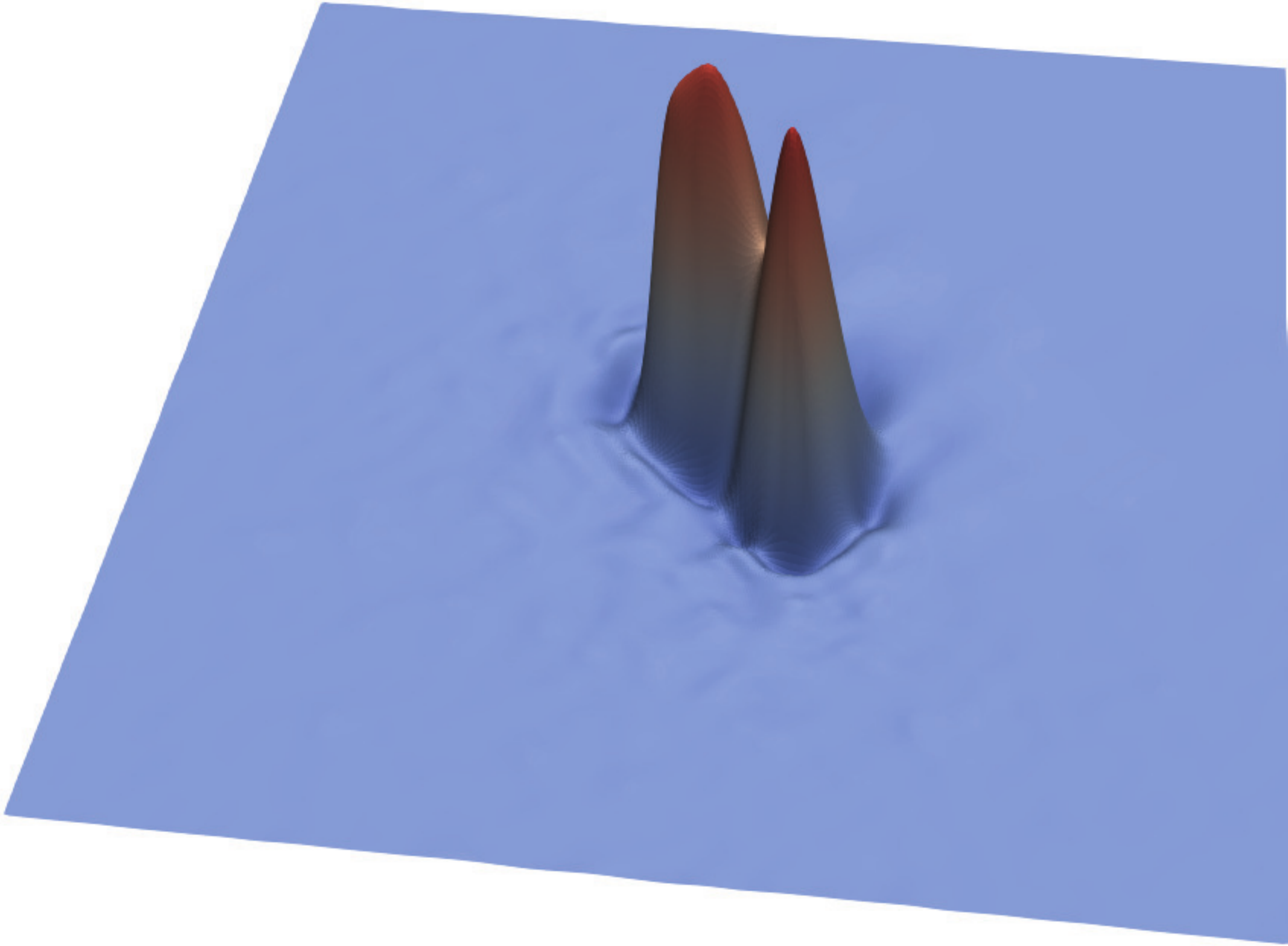}
\caption{B-spline model of the modified polysinc data set with six different values of the regularization threshold $s^*$. Top row, left to right: $s^* = 0.5$ $s^*= 1.0$, and $s^* = 2.0$. Bottom row, left to right: $s^* = 4.0$, $s^* = 8.0$, and $s^* = 16.0$.}
\label{fig:thresh-test-signal}
\end{figure*}

\begin{figure*}
  \centering
  \includegraphics[width=.16\textwidth]{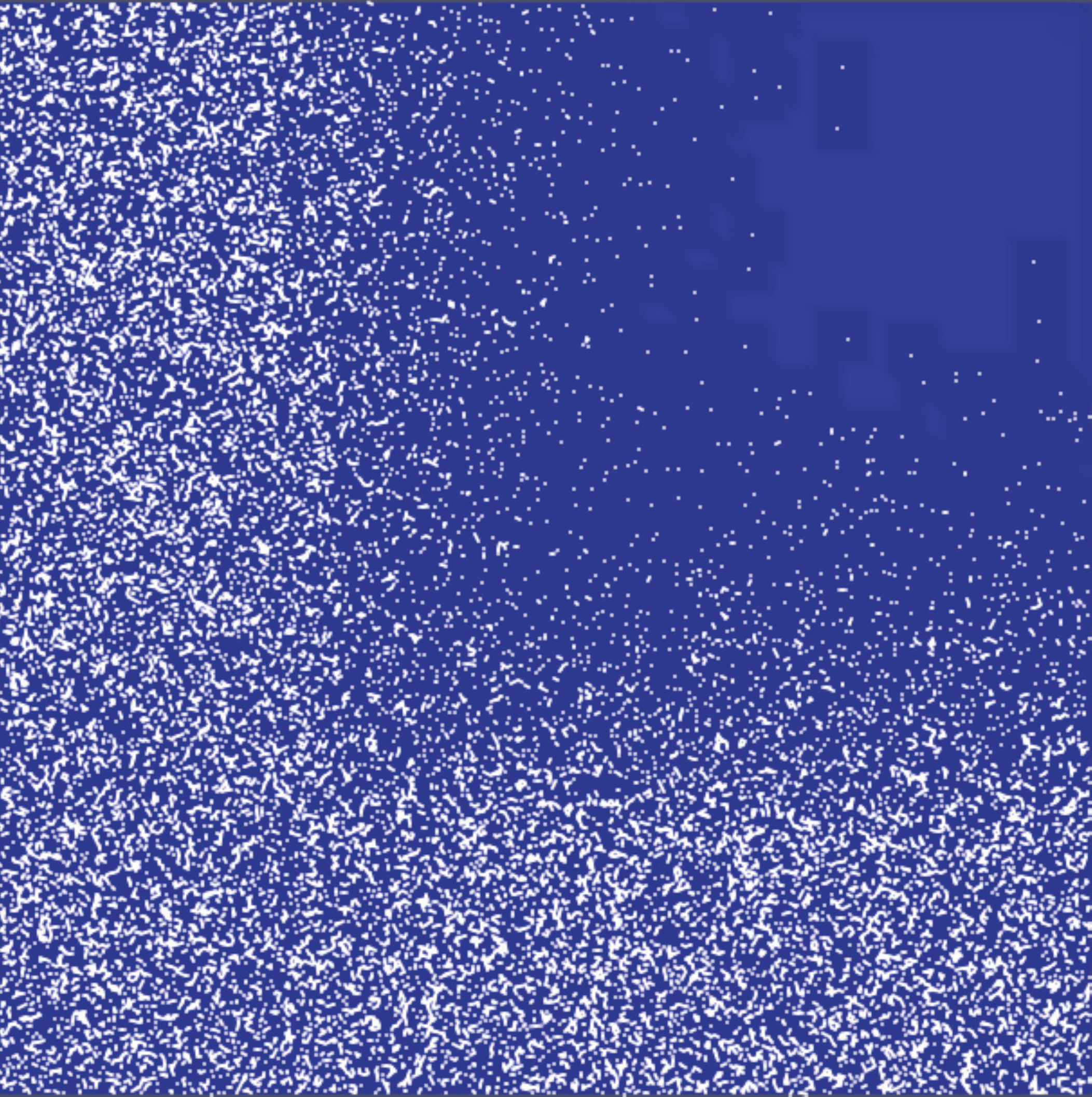}
  \includegraphics[width=.16\textwidth]{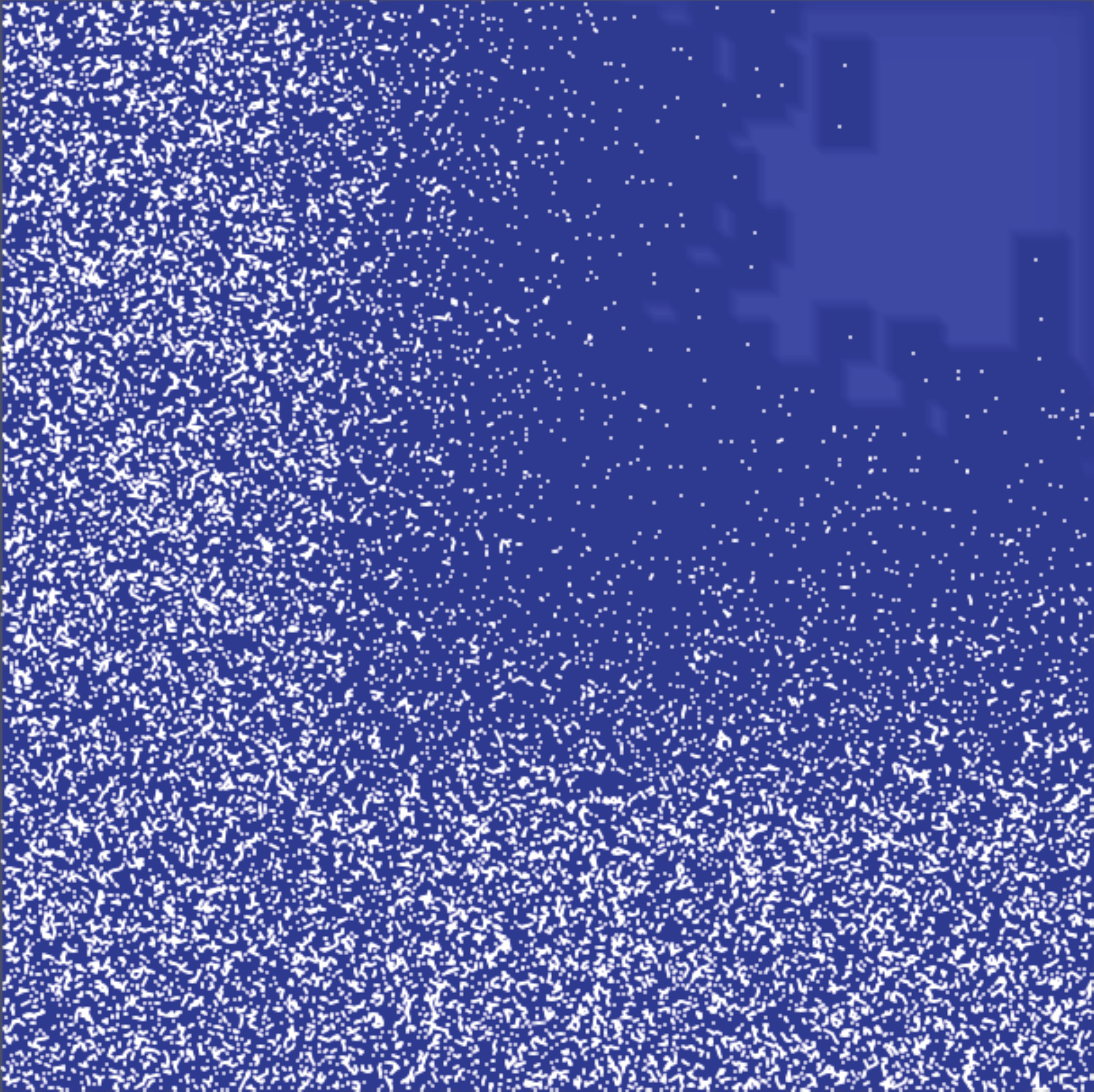}
  \includegraphics[width=.16\textwidth]{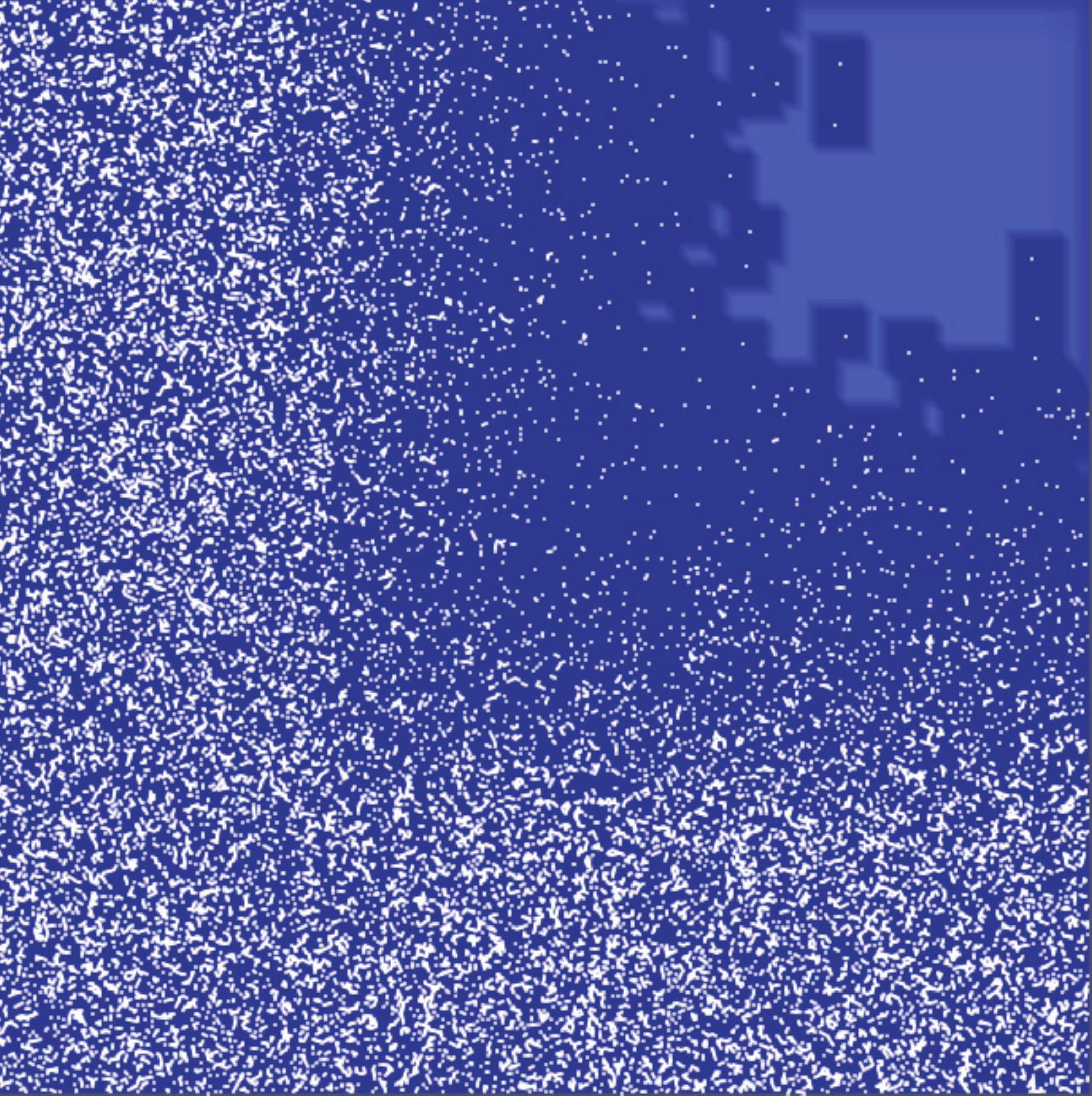}
  \includegraphics[width=.16\textwidth]{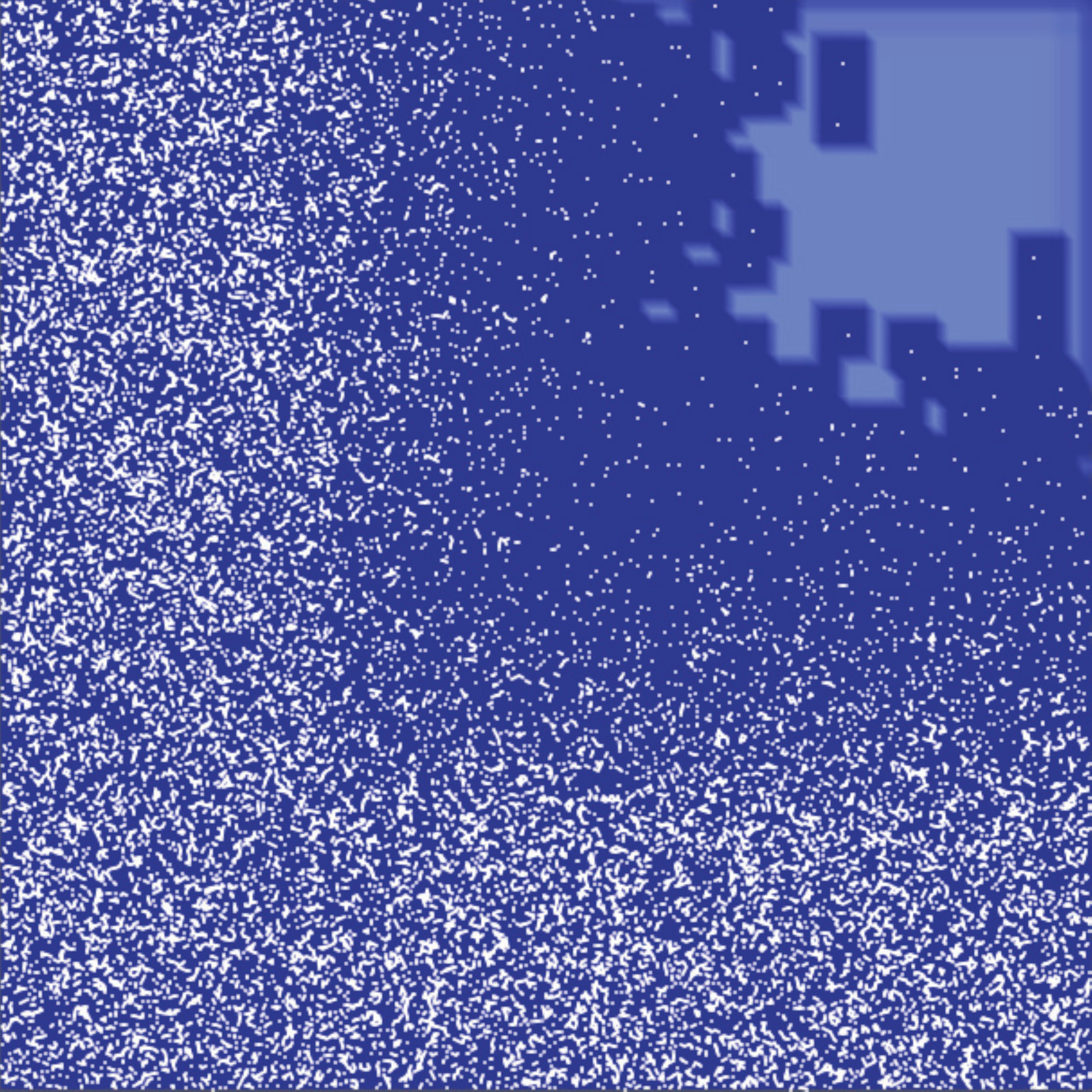}
  \includegraphics[width=.16\textwidth]{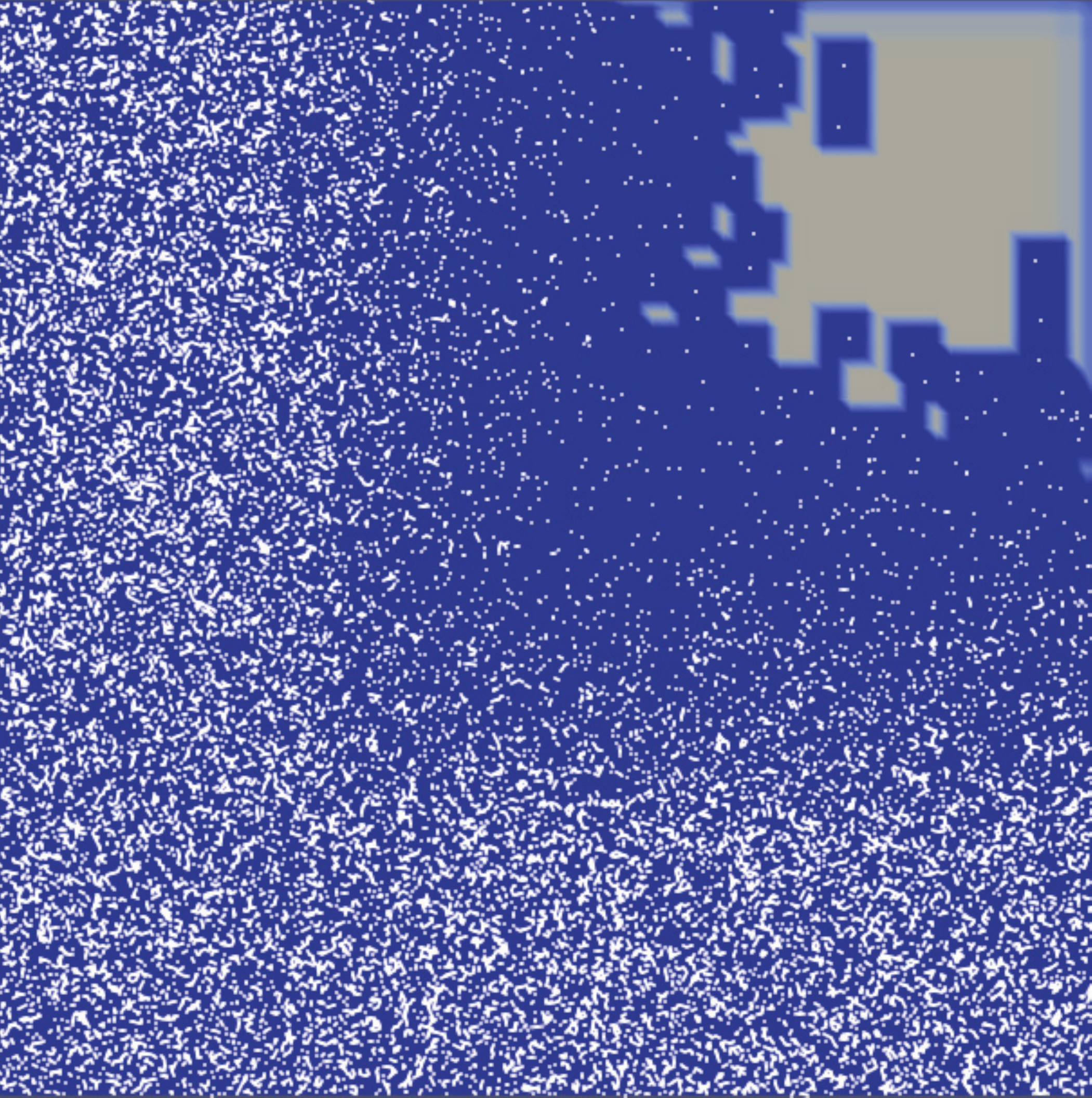}
  \includegraphics[width=.16\textwidth]{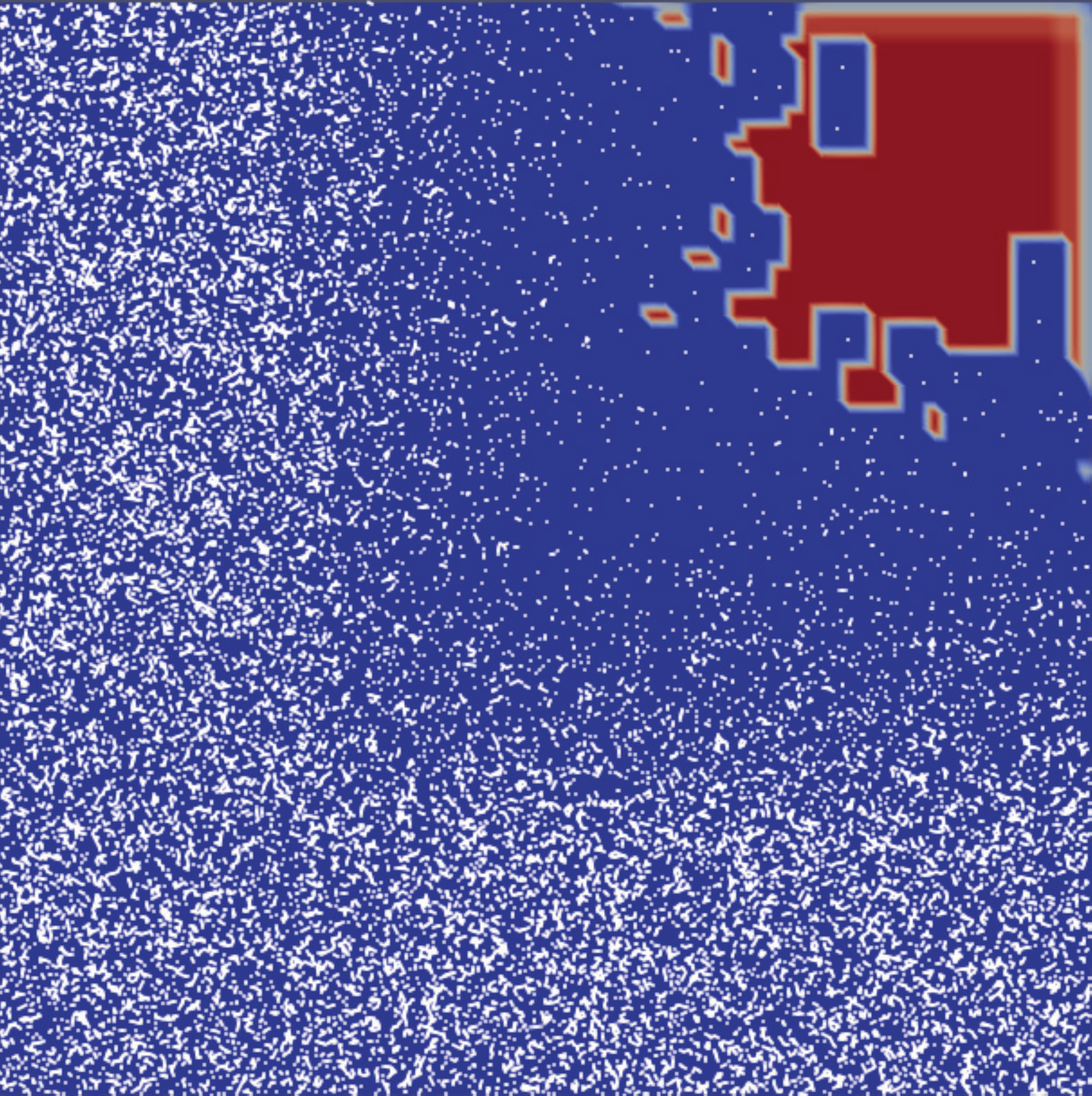}
  \caption{First derivative local regularization strength for the B-spline models in Figure~\ref{fig:thresh-test-signal} of the modified polysinc data set with six different values of the regularization threshold $s^*$. Point locations are overlaid in white. Top row, left to right: $s^* = 0.5$ $s^*= 1.0$, and $s^* = 2.0$. Bottom row, left to right: $s^* = 4.0$, $s^* = 8.0$, and $s^* = 16.0$.}
  \label{fig:thresh-test-strengths-1}
\end{figure*}

\begin{figure*}
  \centering
  \includegraphics[width=.16\textwidth]{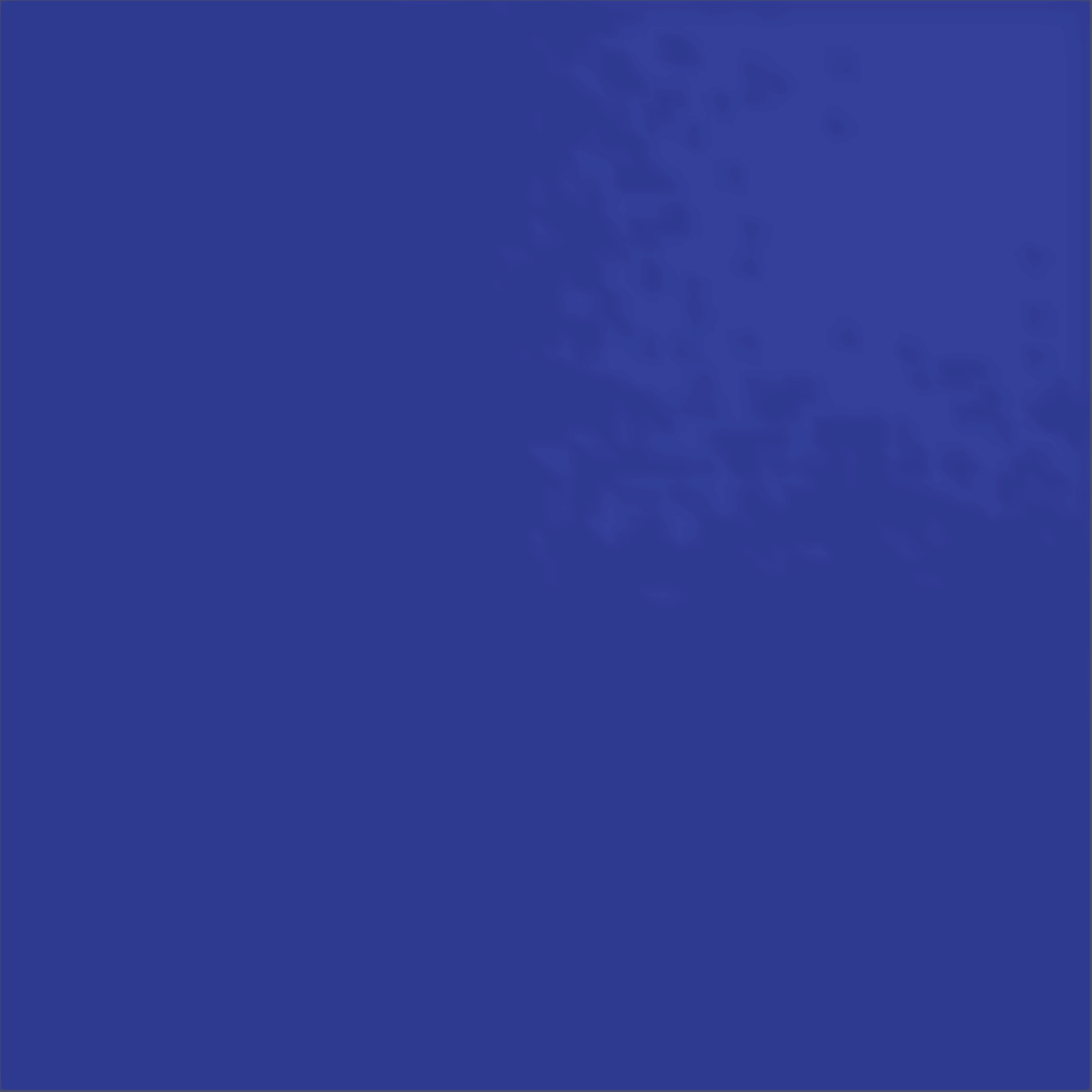}
  \includegraphics[width=.16\textwidth]{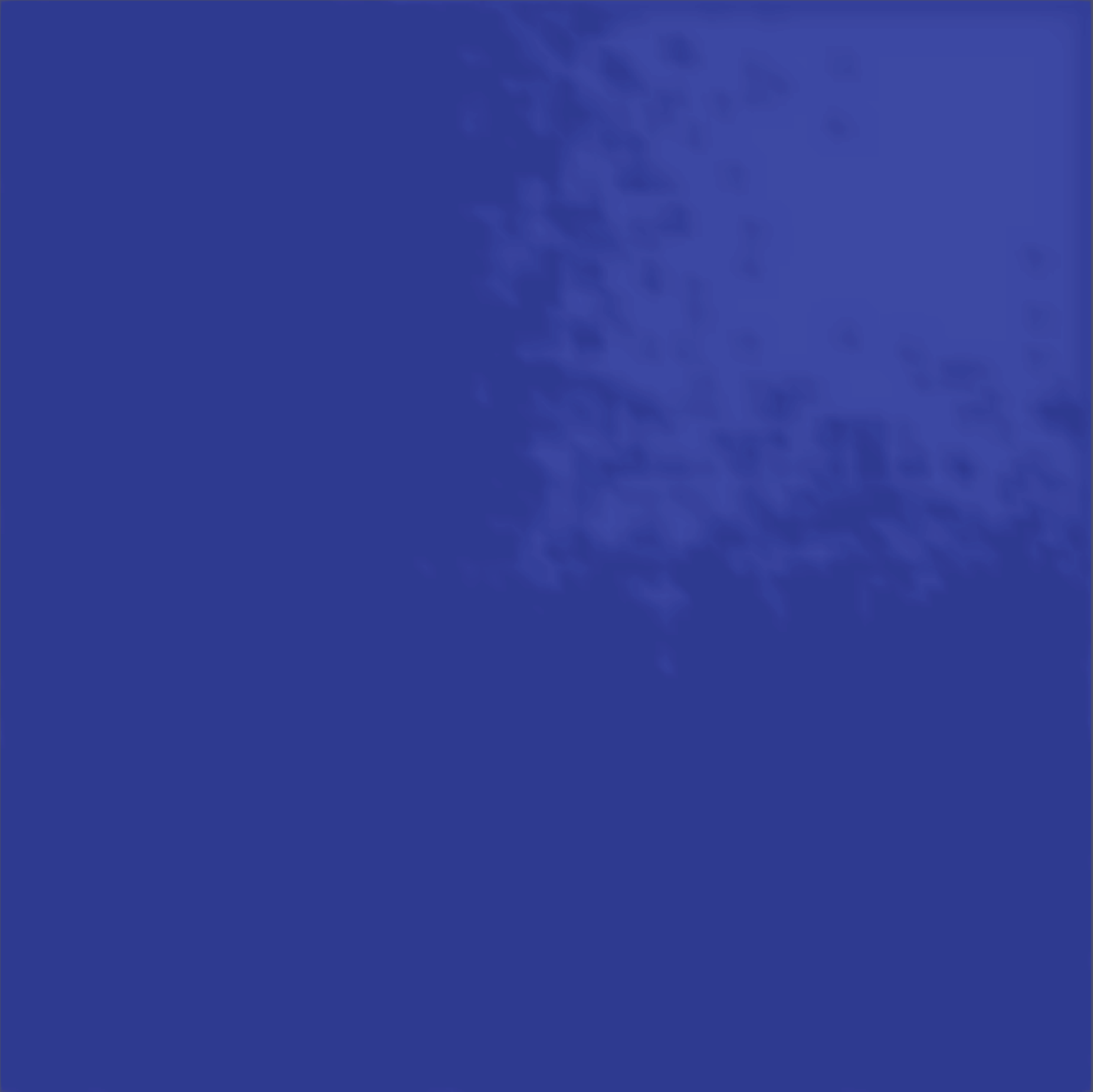}
  \includegraphics[width=.16\textwidth]{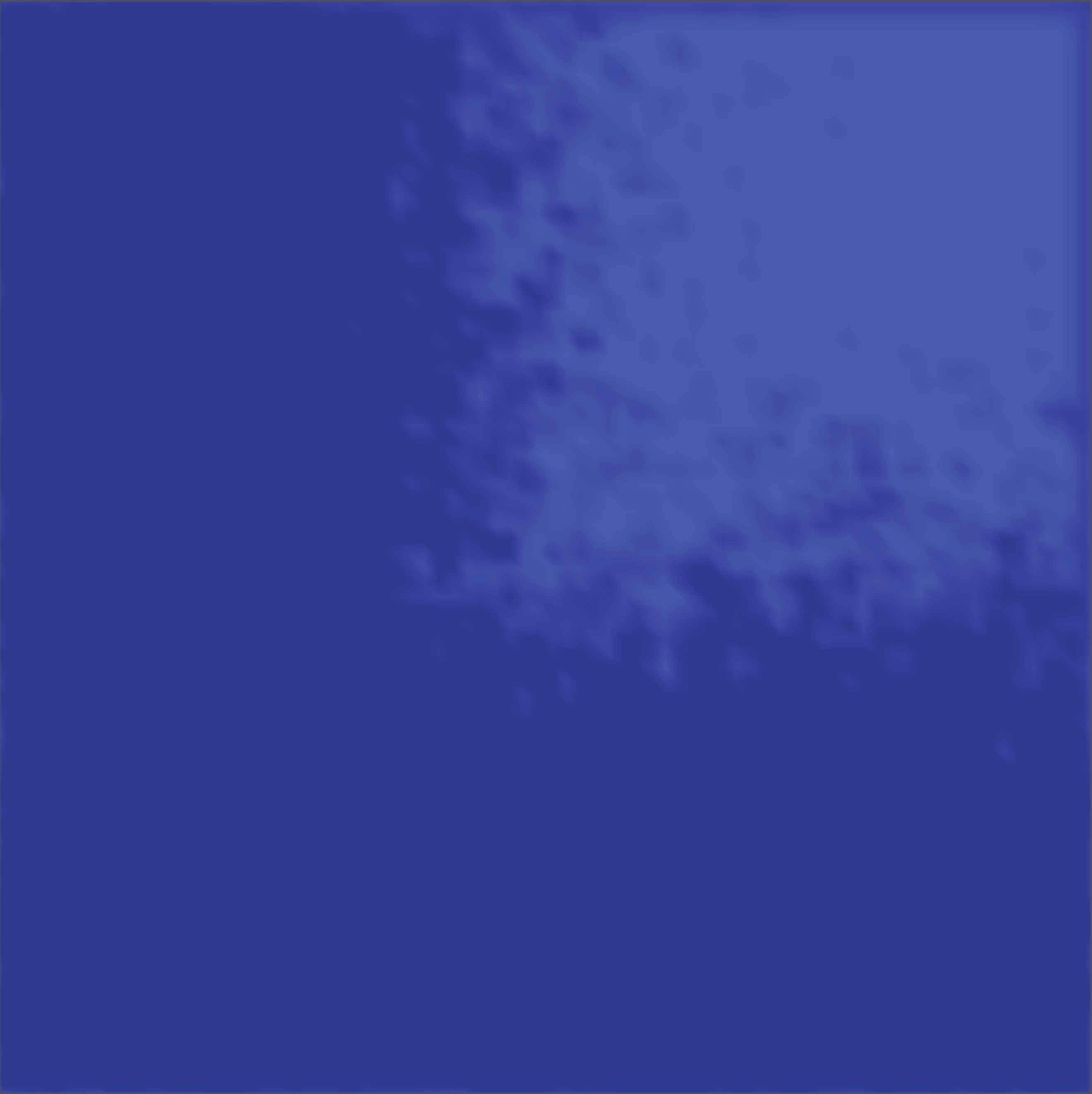}
  \includegraphics[width=.16\textwidth]{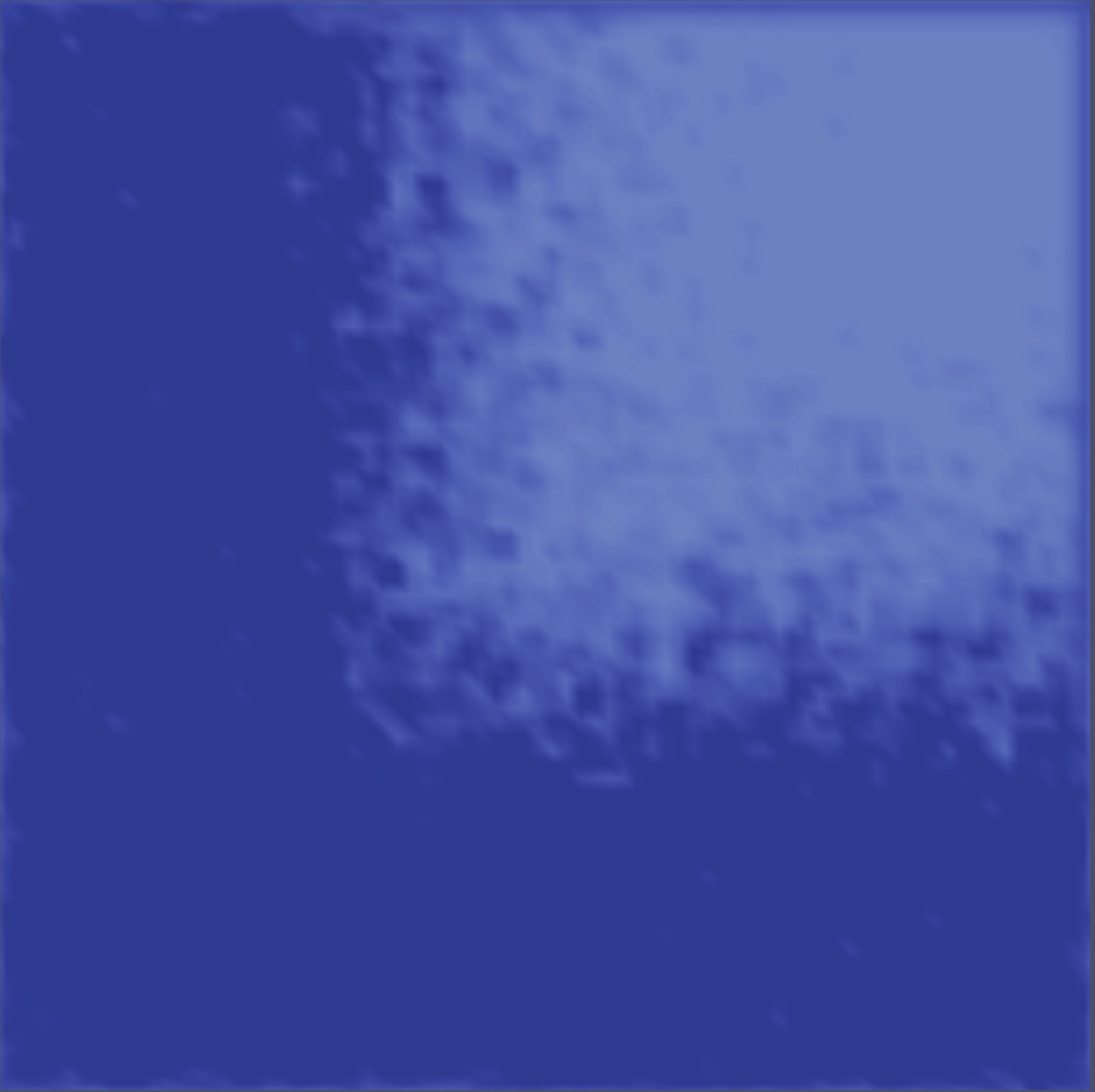}
  \includegraphics[width=.16\textwidth]{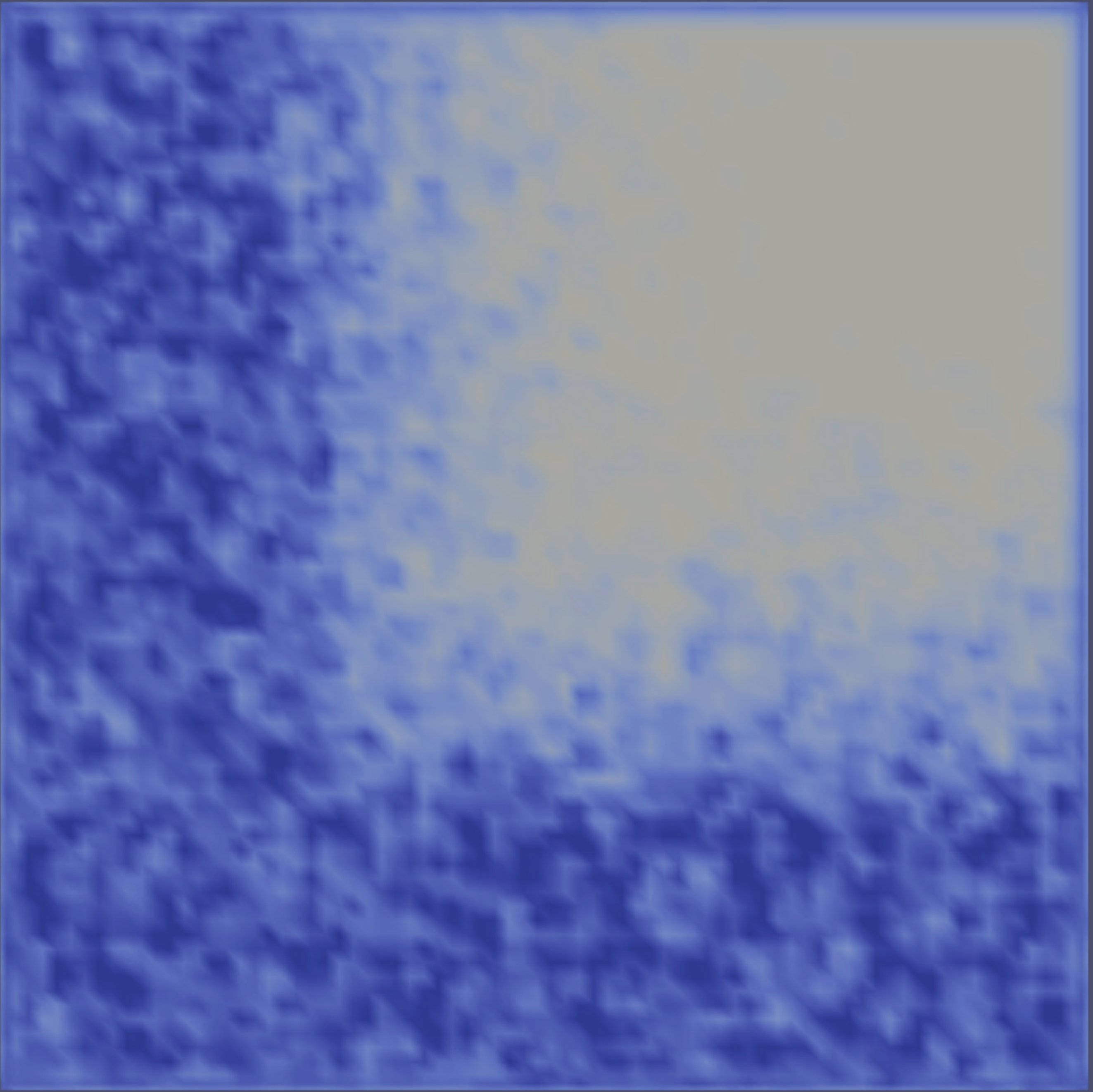}
  \includegraphics[width=.16\textwidth]{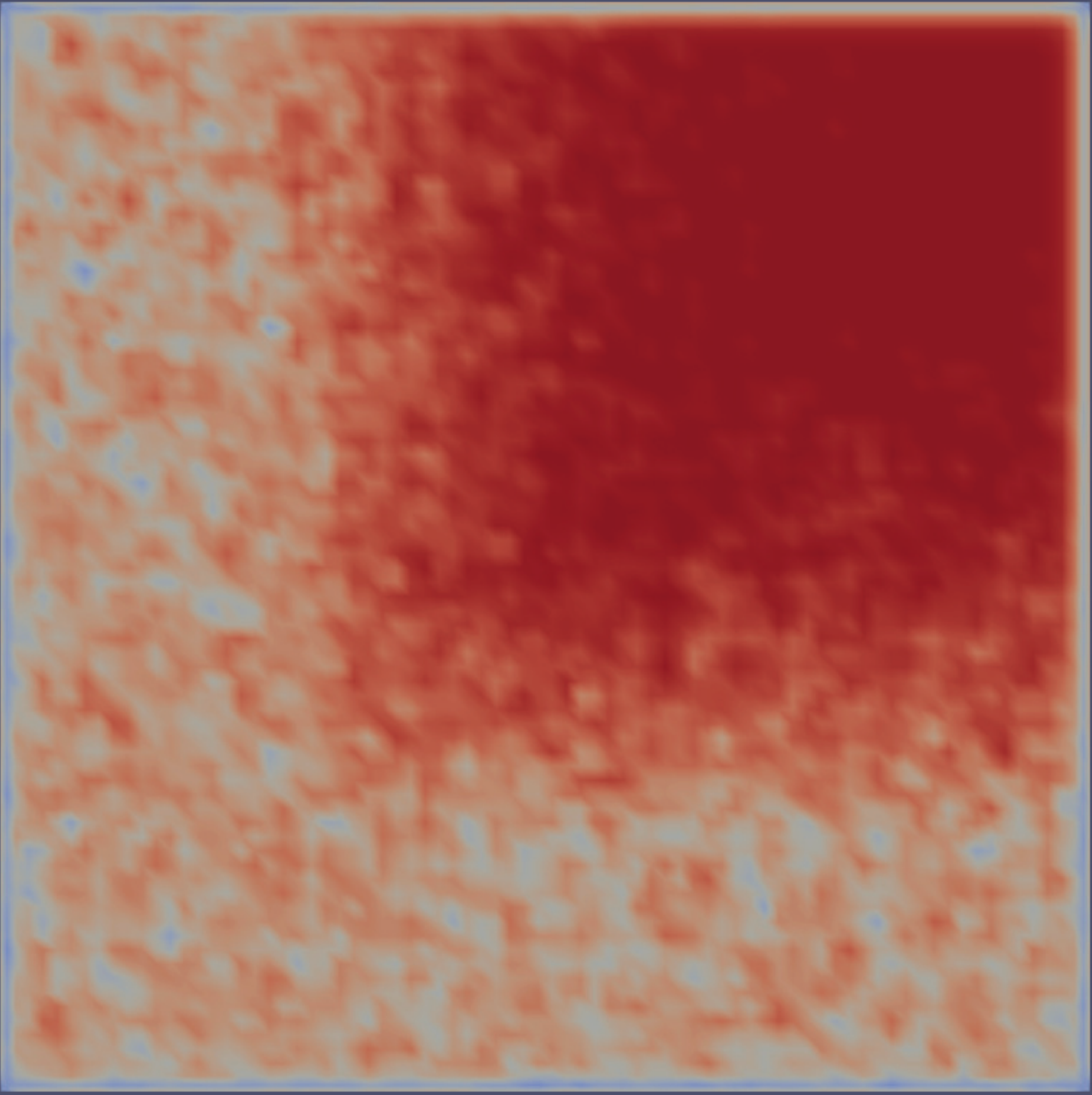}
  \caption{Second derivative local regularization strength for the B-spline models in Figure~\ref{fig:thresh-test-signal} of the modified polysinc data set with six different values of the regularization threshold $s^*$. Top row, left to right: $s^* = 0.5$ $s^*= 1.0$, and $s^* = 2.0$. Bottom row, left to right: $s^* = 4.0$, $s^* = 8.0$, and $s^* = 16.0$.}
  \label{fig:thresh-test-strengths-2}
  \end{figure*}

\subsection{Effect of Regularization Threshold on Model Quality}
A fundamental parameter in the definition of the adaptive regularization method is the regularization threshold, $s^*$. This parameter determines which regions of a domain will be regularized and also sets a maxium regularization strength for the problem. Depending on the requirements of an application and characteristics of a data set, what constitutes the ``best'' value of $s^*$ will vary. However, in all of the examples considered in this article, we found that it was a value of $s^* \in [0,10]$ produced the most faithful model with minimal numerical artifacts. 

To detail the effects of the regularization threshold on model quality, we studied the accuracy and local regularization strengths on a model of the modified polysinc data set over six different values of $s^*$. The complete adaptive regularization method with first and second derivative terms was used to create the model.  The B-spline model was degree 3 with 6,400 total control points ($80\times 80$ grid).

Figure~\ref{fig:thresh-test-signal} shows the result of the adaptive regularization procedure with $s^* \in \{0.5, 1.0, 2.0, 4.0, 8.0, 16.0\}$. For a reference image of the true function, we refer the reader back to Figure~\ref{fig:psinc-2d}. We observe that a regularization threshold of $s^*=0.5$ is clearly inadequate, as the model shows large spurious peaks and oscillates out of the visible frame.  In this case, the model is under constrained. Conversely, the models corresponding to $s^*=8$ and $s^*=16$ may be considerd over-smoothed, with distinctive features of the polysinc function blurred out even in regions with high point density (c.f. Figure~\ref{fig:polysinc-points} for the point distribution). In this example, values of $s^*$ between 1 and 4 produced the best trade off between accuracy and well-conditioning.

Table~\ref{tab:thresh-test} displays the analytical errors corresponding to each of the six regularized models. To compute theses error metrics, we sampled the B-spline model on a high-resolution regular grid ($400\times 400$) and compared the model's output with the true function value given by the polysinc function. Again, we see the worst $L^\infty$ errors occur at the extremes.  The $L^2$ errors do not change significantly due to the fact that all six models fail to adequately capture the low-amplitude, high-frequency oscillations at the edges of the domain. Because poorly-constrained models tend to produce highly localized spurious peaks, the $L^\infty$ error  generally serves as a better indicator of quality than the $L^2$ error, which is less sensitive to high residuals in small neighborhoods.  The difference in behvior seen in Table~\ref{tab:thresh-test} aligns with this observation.

Finally, Figure~\ref{fig:thresh-test-strengths-1} and Figure~\ref{fig:thresh-test-strengths-2} contain plots of the local regularization strengths as a function of position. Figure~\ref{fig:thresh-test-strengths-1} shows the first derivative local regularization strength, which is zero except in the top-right corner, no matter the value of $s^*$. This is what we expect from our method, since $\mathbf{\Lambda_1}$ is defined to vanish wherever there is at least one input point nearby. Here, increasing $s^*$ simply gradually increases the local regularization strength in a fixed region. In contrast, the second derivative local regularization strength (shown in Figure~\ref{fig:thresh-test-strengths-1}) is near zero when $s^*$ is small, but gradually spreads throughout the domain as $s^*$ increases. For this reason, we can control \emph{how much} of the domain is regularized by increasing $s^*$, with the regions with lowest point density being regularized ``first.'' Finally, we remark here that the accuracy deficiencies seen in Figure~\ref{fig:thresh-test-signal} can be directly explained by looking at the second derivative local regularization strength. When $s^* = 0.5$ and the model exhibits spurious peaks, we see that those peaks occur in regions of low point density where the local regularization strength is zero. In contrast, when $s^*=8$ and $s^*=16$ and the model is oversmoothed, we observe that the local regularization strength is very high even in the lower left corner where point density is highest and regularization is unnecessary.

\begin{table}
  \caption{Analytical errors as a function of the regularization threshold, $s^*$ for the models shown in Figure~\ref{fig:thresh-test-signal}. Errors were computed by comparing the modeled value to the true function value on a $400\times 400$ regular grid of points. $L^2$ Error is also known as the root mean-squared error, $L^\infty$ error as the maximum pointwise error.}
  \centering
  \begin{tabularx}{\linewidth}{@{}l @{\hspace{5pt}}|@{\hspace{3pt}} 
    *6{>{\centering\arraybackslash}X}@{}}
    $s^*$ & 0.5 & 1.0 & 2.0 & 4.0 & 8.0 & 16.0\\ \hline
    $L^2$ Error & 0.639 & 0.265 & 0.246 & 0.267 & 0.344 & 0.438\\
    $L^\infty$ Error & 22.3 & 6.01 & 3.93 & 3.40 & 3.88 & 5.12\\
  \end{tabularx}
  \label{tab:thresh-test}
\end{table}

\subsection{Extrapolation into Unconstrained Regions}
When a large region of the domain does not contain any data points to constrain the best-fit B-spline problem, the least-squares minimization will be ill-posed and the resulting model can exhibit extreme oscillations. However, data sets with empty regions or ``holes'' are very common in scientific and industrial applications. For example, some climate models measure ocean temperatures or land temperatures, but not both simultaneously. Industrial simulations often model objects with irregular boundaries, and data from physics simulations are shaped by the locations of detectors.

Although empty regions are usually omitted in subsequent analysis, it is still important to understand and control the behavior of a B-spline model in empty regions. Extreme oscillations near the boundary of a hole can distort the derivative of the model away from this boundary. In addition, attempting to compute simple statistics about the model (minimum, maximum, mean) can be biased if the model exhibits unpredictable behavior in empty regions.

\begin{figure}
  \centering
  \includegraphics[width=\linewidth]{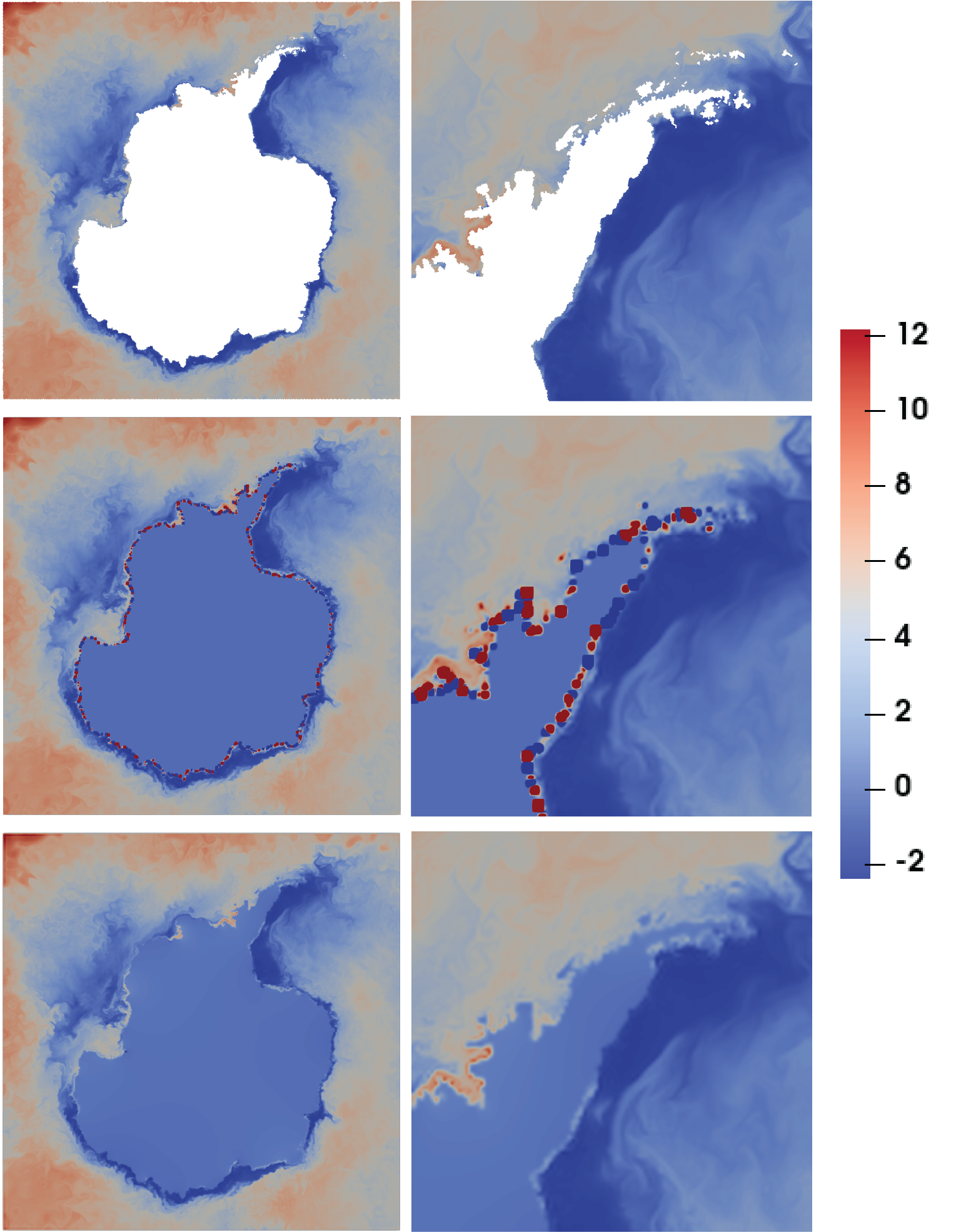}
  \caption{Ocean temperature simulation around Antarctica. Top row: Input data on a hexagonal mesh (left), B-spline model without regularization (center), B-spline model with adaptive regularization (right). Bottom row: Detail view from the top row.} 
  \label{fig:antarctica}
\end{figure}

Figure~\ref{fig:antarctica} shows a scenario from a simulation of ocean temperatures. The data set is centered around the continent of Antarctica, for which no temperature values are given. Exhibited in the figure are two models, both of degree 2 with a 400 $\times$ 400 grid of control points. The unregularized model at center oscillates between $\pm 10^8$ along the coast (while the input data range from -2 to 10). In contrast, the regularized model (with $s^* = 5$) smoothly transitions at the coast to a near-constant value over the landmass. Regularization was performed with constraints on first and second derivatives, as described at the end of Section~\ref{sec:method}.

\begin{figure}
  \centering
  \includegraphics[width=100pt, height=100pt]{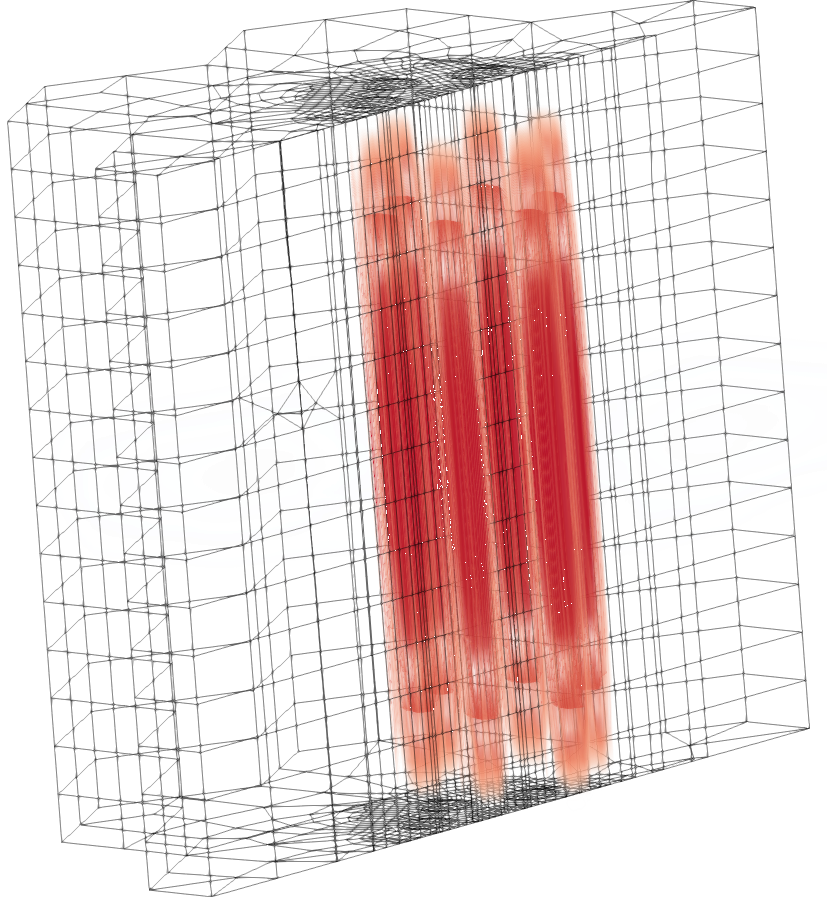}
  \includegraphics[width=100pt, height=100pt]{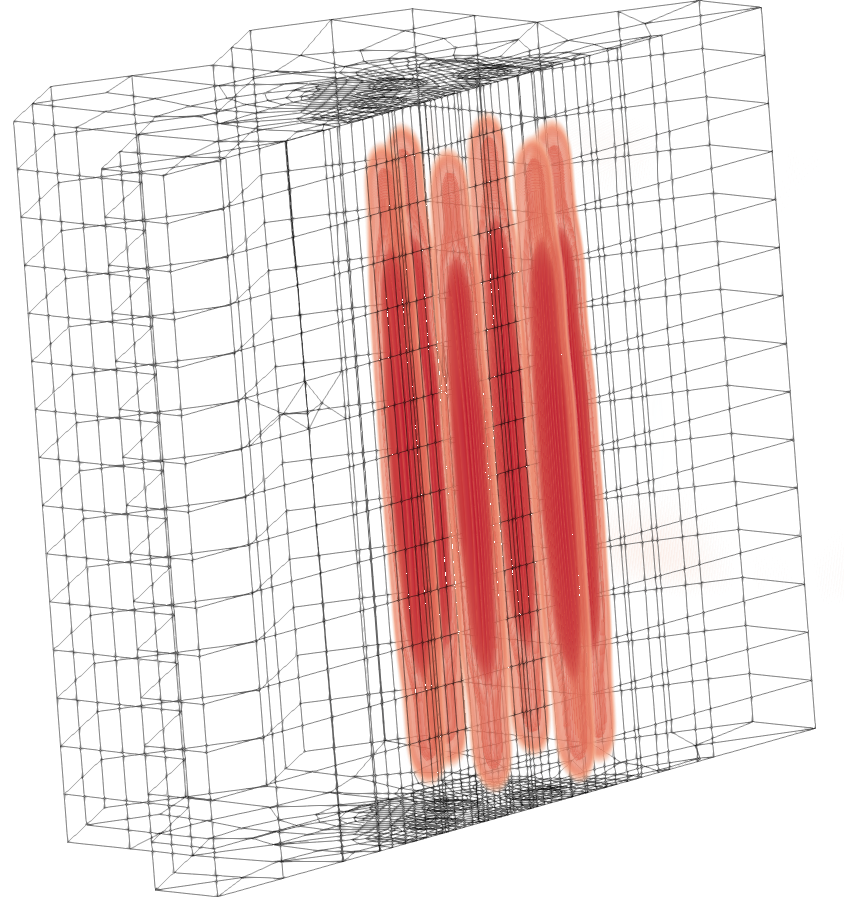}
  \includegraphics[width=100pt, height=100pt]{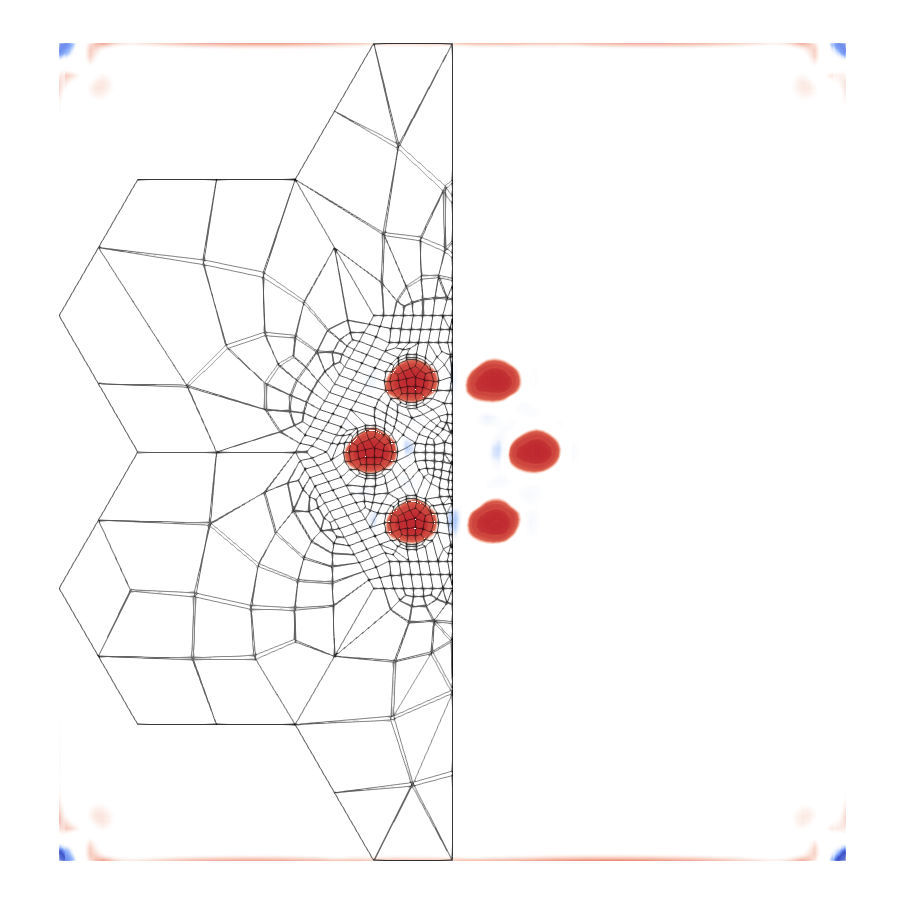}
  \raisebox{10pt}{\includegraphics[width=35pt]{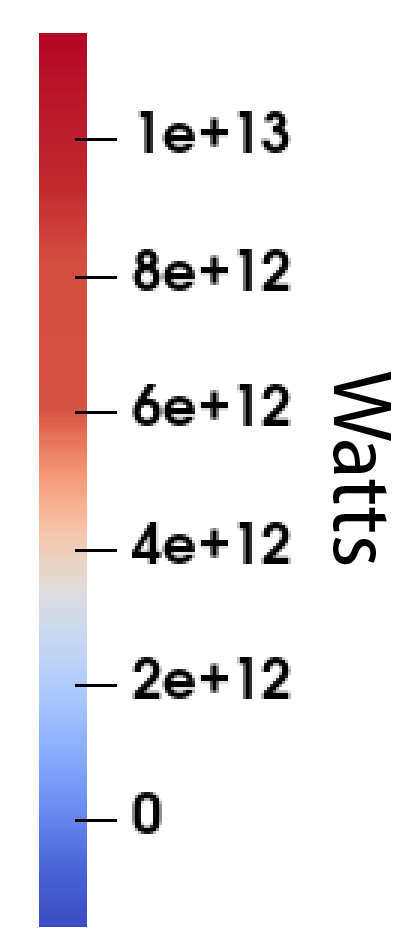}}
  \caption{Power production (Watts) in a nuclear reactor simulation. Input data (left), spline model (center), spline model top view (right).}
  \label{fig:nuclear}
\end{figure}

We next consider a three-dimensional data set representing power produced in a component of a nuclear reactor (Figure~\ref{fig:nuclear}). The data are contained inside a hexagonal prism, but the B-spline model is defined on the bounding box of this prism. Hence the corners of this box are devoid of any data. 
Without regularization, the least-squares minimization does not converge, so we exhibit only our regularized model (with $s^*=10$) in Figure~\ref{fig:nuclear}. The adaptive regularization method (Figure~\ref{fig:nuclear}, center and right) produces an accurate model of the six interior ``pins'' and is well defined in the corner regions. Some artifacts are observed in the corners, but they are not significant enough to affect the interior of the model. Without regularization, the condition number of the system is infinite; with adaptive regularization, the condition number is $2.18\times10^5$. The model was produced by constraining first and second derivatives simultaneously.

\section{Future Work}\label{sec:future-work}
In the construction of our method, we imposed constraints on the first and second derivatives of the B-spline in regions where the density of input data was low or vanishing. However, these additional constraints could have been based on different orders of derivative or been unrelated to derivatives altogether. For instance, a different type of constraint would be one that penalizes deviation from a known baseline value. Another option would be to penalize B-spline values that exceed a given range (for example, the original bounds of the input data). We remark that our procedure for adaptive regularization of tensor product B-splines is separate from the type of artificial constraint imposed. Depending on the application, different constraints may be more useful, and our method allows for those to be used instead.

Another direction for future research is to further adapt the regularization strength based on the local properties of the input data. At present, the same regularization threshold is used throughout the domain. This definition makes the method very easy to use, but more flexible or complex schemes are possible. For instance, different regularization thresholds could be used separately for first and second derivative regularization, potentially allowing for a better balancing of these two kinds of constraints. In addition, the method as presented here varies the regularization strength from $0$ to the value $s^*/\widetilde{s}_j$, as defined in Equation~\eqref{eq:lambda-j}. In applications where greater user interaction is expected, the method may be adapted to enforce a customized local minimum or maximum level of smoothing in different subregions of the domain. Such flexibility is especially useful when the data (or subregions of the data) are known to be noisy or, conversely, highly-sensitive to smoothing.

\section{Conclusions}\label{sec:conclusions}
Modeling unstructured data sets with tensor product B-splines can be difficult due to the ill-conditioning of the fitting problem. In general, data sets with large variations in point density or regions without data exacerbate this problem to the point that artificial smoothing is necessary. However, smoothing an entire model can wash out sharp features in the data.

We introduced a regularization procedure for B-spline models that preserves features by adapting the regularization strength throughout the domain. Our method automatically varies the smoothing intensity as a function of input point density and relies on a single user-specified parameter, which we call the regularization threshold. We observe that adaptive regularization performs better than typical uniform regularization schemes that may over-smooth some regions while under-smoothing others. We also showed that our method can fit B-spline models to data sets with regions of extremely sparse point density and remain well-defined even in areas without data points. Overall, adaptive regularization of B-spline models produces smooth and accurate models for data sets which would otherwise be difficult to fit.

%
%
\bibliographystyle{elsarticle-num-names} 
\bibliography{mfa-reg}
\end{document}